\documentclass[draftcls,onecolumn,12pt]{IEEEtran}

\usepackage{amssymb,times}
\usepackage{cite, epsfig}
\usepackage{verbatim, amsfonts,graphicx,enumerate}     
\usepackage{bm,array}
\usepackage[cmex10]{amsmath}
\usepackage{graphicx}				
\usepackage{amsthm}
\usepackage{epstopdf}
\usepackage{float}
\usepackage{algorithm,algorithmic}
\usepackage{framed}
\usepackage[font=small]{caption}
\usepackage{subcaption}
\usepackage{color,url}
\usepackage{multirow}
\usepackage[normalem]{ulem}
\theoremstyle{plain}
\newtheorem{thm}{Theorem}
\newtheorem{lem}{Lemma}
\newtheorem{cor}{Corollary}

\theoremstyle{definition}

\theoremstyle{remark}

\DeclareMathOperator*{\argmax}{arg\,max}
\def \x {{\mathbf{x}}}
\def \Q {{\mathbf{Q}}}

\def \p {{\mathbf{p}}}
\def \s {{\mathbf{s}}}
\def \f {{\mathbf{f}}}

\def \ft {{\mathbf{f}_t}}
\def \f {{\mathbf{f}}}
\def \fa {{\mathbf{f}}}

\def \st {{\text{s. t. }}}
\def \Xc {{\mathcal{X}}}
\def \Pc {{\mathcal{P}}}

\def \lam {{\bm{\lambda}}}
\def \lamt {{\bm{\lambda}_t}}
\def \ls {{\bm{\lambda}^\star}}
\def \lt {{\tilde{\bm{\lambda}}}}

\def \e {{\mathbf{e}}}

\def \u {{\mathbf{u}}}
\def \v {{\mathbf{v}}}

\def \Rn {{\mathbb{R}}}
\def \Nn {{\mathbb{N}}}

\def \O {\mathcal{O}}
\def \bg {{\boldsymbol{\varphi}}}

\def \E {{\mathbf{E}}}				
			
\def \mu {\epsilon}
\def \st {{\mathbf{s}_t}}
\def \s {\mathbf{s}}

\def \T {{\mathcal{T}}}
\def \eps {{\epsilon}}
\def \pa {{(\mathbf{P}_1)}}

\providecommand{\Ex}[1]{\mathbb{E}\left[#1\right]}
\providecommand{\Eh}[1]{\mathbb{E}_h\left[#1\right]}
\newcommand{\colb}[1]{\textcolor{black}{#1}}

\providecommand{\pl}[1]{\mathcal{P}_{\Lambda}\left(#1\right)}
\providecommand{\pll}[1]{\mathcal{P}_{\Lambda}\left(#1\right)}
\providecommand{\ip}[1]{\langle#1\rangle}

\providecommand{\abs}[1]{\left|#1\right|}
\providecommand{\norm}[1]{\left\|#1\right\|}

\allowdisplaybreaks
\begin{document}
	\title{\vspace{-13mm}\colb{Network Resource Allocation via Stochastic Subgradient Descent: Convergence Rate}	}
	\author{Amrit S. Bedi,~\IEEEmembership{Student Member, IEEE} and Ketan Rajawat,~\IEEEmembership{Member, IEEE} \thanks{Manuscript submitted \today. This work was supported by the Indo-French Centre for the Promotion of Advanced Research-CEFIPRA. The authors are with the Department of Electrical Engineering, IIT Kanpur, Kanpur (UP), India 208016 (email: \texttt{\{amritbd,ketan\}@iitk.ac.in}).}\vspace{-18mm}}
	\vspace{0mm}
	\maketitle

	\begin{abstract}
		This paper considers a general stochastic resource allocation problem that arises widely in wireless networks, cognitive radio, networks, smart-grid communications, and cross-layer design. The problem formulation involves expectations with respect to a collection of random variables with unknown distributions, representing exogenous quantities such as channel gain, user density, or spectrum occupancy. We consider the constant step-size stochastic dual subgradient descent (SDSD) method that has been widely used for online resource allocation in networks. The problem is solved in dual domain which results in a primal resource allocation subproblem at each time instant. The goal here is to characterize the non-asymptotic behavior of such stochastic resource allocations in an almost sure sense.
		 It is well known that with a step size of $\epsilon$, {SDSD} converges to an $\mathcal{O}(\epsilon)$-sized neighborhood of the optimum. In practice however, there exists a trade-off between the rate of convergence and the choice of $\epsilon$. This paper establishes a convergence rate result for the SDSD algorithm that precisely characterizes this trade-off. {Towards this end, a novel stochastic bound on the gap between the objective function and the optimum is developed. The asymptotic behavior of the stochastic term is characterized in an almost sure sense, thereby generalizing the existing results for the {stochastic subgradient} methods.} For the stochastic resource allocation problem at hand, the result explicates the rate with which the allocated resources become near-optimal. As an application, the power and user-allocation problem in device-to-device networks is formulated and solved using the {SDSD} algorithm. Further intuition on the rate results is obtained from the verification of the regularity conditions and accompanying simulation results. 
	\end{abstract}
		\vspace{-5mm}
	\begin{IEEEkeywords}
		\vspace{-5mm}
		Stochastic subgradient, constant step-size, stochastic resource allocation, D2D communication.
	\end{IEEEkeywords}
	\vspace{-8mm}
	\section{Introduction}
	\IEEEPARstart{R}{esource} allocation is a fundamental problem in economic theory that finds application in the design of wireless communication protocols \cite{resource_alloca_tut_1_neely}, smart grid systems \cite{alloc_smart_tut_1}, and scheduling algorithms \cite{alloc_sched_tut_1}. From an optimization perspective, the goal is to find the optimal allocation variables such as transmit power, bandwidth, operational schedule, or facility locations, so as to maximize the user satisfaction, minimize the cost, and satisfy all system constraints. The \emph{stochastic} resource allocation problem arises in scenarios where the optimization problem includes random parameters with unknown distributions \cite{stoch_res_prob_3}. For such problems, the goal is to find an allocation \emph{policy} that is feasible and optimal, \colb{on average} \cite{separation} or with high probability \cite{marques3}.  Since the policy variable may be infinite dimensional , the problem is more tractable in the dual domain {due to finite number of constraints}, an aspect exploited by a number of algorithms; see e.g., \cite{crolayketan,ale10,ofdm} and references therein.
	
	This paper focuses on the so-called online algorithms, where the allocation must occur every time the random parameter is realized and revealed. For each realization, the resource allocation adheres to the operational or ``box'' constraints, while the overall allocation policy is only asymptotically feasible and optimal. The dual problem may then be solved using the stochastic subgradient descent method, whose asymptotic behavior is well-known \cite{nedic2001convergence,bertsekas2011incremental}. Further justification for solving the problem in the dual domain was provided in \cite{separation,crolayketan}, where it was shown that such stochastic problems have zero duality gap under some mild conditions. The asymptotic feasibility and optimality of the allocated resources via primal averaging was also established in \cite{ale10}. In a similar vein, the relationship between the stochastic and 'averaged' dual iterates for the power and subcarrier allocation problem in OFDM was established in \cite{ofdm}. 
	
	In resource allocation problems, it is possible for the environmental variables to change abruptly. This motivates the use of constant step sizes in stochastic algorithms, that obviate the need to restart the iterations whenever such a change occurs \cite{Luo_LMS_rate}. With a constant step size of $\eps$ however, it is well known that the stochastic iterates converge only to an $O(\eps)$-sized neighborhood of the optimum \cite{ale10}. On the other hand, making $\eps$ arbitrarily small is also impractical, since it results in a slow convergence rate \cite{bottou1998online,bottou2012stochastic}. The aforementioned trade-off between the rate of convergence and the asymptotic performance of the constant step-size \colb{stochastic dual subgradient descent} {(SDSD)} algorithm is an important aspect that has not been studied explicitly.
	
The goal of this paper is to rigorously characterize the convergence rate of the {SDSD} algorithm in an almost sure sense. 	{
		The key contribution of the paper is the development of stochastic bounds on the iterates produced by {SDSD} method,
		that explicate the role played by $\eps$ in ‘forgetting’ the initial conditions, and coming close to the optimum.} To this end, the iterations are divided into epochs of duration $1/\epsilon$, and the optimality gap is analyzed for both fixed and arbitrarily small $\epsilon$. The main result of the paper is that the stochastic component of this gap goes to zero almost surely, either as the number of epochs go to infinity with fixed $\epsilon > 0$, or as $\epsilon$ itself goes to zero. The bounds developed here specialize to the known asymptotic results, and are generally applicable to any problem solved via the {SDSD} method. {To the best of our knowledge, these are also the first such convergence rate results for the constant step size {SDSD} algorithm.} Corresponding results for the diminishing step size stochastic subgradient exist, but cannot be readily extended to the present case \cite{bach2011non}. Likewise, the analysis in \cite{nedic2009approximate} can be extended to yield rate results that hold \colb{on average}, but does not yield almost sure bounds. The analysis in the present work makes use of the strong law of large numbers directly, and is completely different from that in \cite{bach2011non,nedic2009approximate}.
	%The bounds consist of deterministic and stochastic terms, each of which decay differently, leading to different interpretations for vanishing and fixed $\eps$. 

	As the second contribution, it is shown that the convergence rate results are readily applicable to the stochastic resource allocation problem of interest here. To further demonstrate the usefulness of the bounds, the paper details a contemporary application that uses mobile caching for improving service via device-to-device (D2D) communication \cite{ji2016wireless,archi_D2D}. To this end, we consider the D2D edge caching framework where willing users offer data connectivity to highly mobile users experiencing spotty coverage. 
The problem is well-motivated for vehicular users who may download data from other users residing near the highway. 
	 The corresponding resource allocation problem is shown to satisfy the required regularity conditions, thereby demonstrating the flexibility afforded by the {SDSD} algorithm.
	
	This paper is organized as follows. Sec. \ref{rel} lists some of the related work in this area, providing context to the current work. \colb{ Sec. \ref{probfor} starts with detailing the D2D edge caching problem and formulates the general network resource allocation problem.} Sec. \ref{solu_dual} discusses the various solution methodologies in the literature, including the {stochastic subgradient descent (SSD)} framework. Sec. \ref{conv} provides the main results of the paper, stating the convergence rate results for both primal and dual problems. Sec. \ref{d2d} further develops the D2D examples introduced in Sec. \ref{probfor}, and verifies the different conditions required for the convergence results to hold. Simulation results on D2D example are provided in Sec. \ref{num} and Sec. \ref{conc} concludes the paper.
	
The notation used here is as follows. Boldface letters denote column vectors, for which the inequalities and equalities are defined component-wise. The set of all real numbers is denoted by $\Rn$, and likewise the sets of non-negative reals, positive reals, and $K$-dimensional real vectors are denoted by $\Rn_{+}$, $\Rn_{++}$, and $\Rn^K$, respectively. Time indices are denoted by the subscripts $t$ and $\tau$. For a vector $\x$, $[\x]_i$ denotes its $i$-th entry, $\norm{\x}$ denotes its $\ell_2$ norm, $\norm{\x}_p$ denotes its $p$-th norm, for $p \in \Rn_{+}$, and $\x^T$ denotes the transposed row vector. The expectation operator is denoted by $\Ex{\cdot}$ and the inner product is denoted by $\ip{\cdot,\cdot}$. Finally, $[c]^a_b = \min(\max\{c,b\},a)$ and $[c]_{+}:=[c]_0^{\infty}$. 
	\vspace{-6mm}
	\subsection{Related work}\label{rel}
	Stochastic approximation algorithms have a long history, going back to the prototypical adaptive filtering algorithms by Robbins and Monro \cite{robbins1951} and Widrow and Stearns \cite{Widrow}, and have been studied extensively in the context of least mean square (LMS) and recursive least-squares (RLS) algorithms \cite{sayed2011adaptive}. Stochastic gradient and subgradient methods have since been applied to neural networks \cite{bottou_91c_neural_learning}, parameter tracking \cite{kushner1994analysis}, large-scale machine learning \cite{bottou2010large,bach2011non}, and resource allocation problems \cite{stoch_res_prob_3}. Convergence of these algorithms is well known for various choices of the step-size parameter \cite{optml}. Convergence rate of the {stochastic subgradient descent} algorithm has been established for diminishing step size rules via {non-asymptotic} analysis \cite{bach2011non}. However, not much is known about the convergence rate of the constant step-size counterpart, except for the fact that it exhibits linear convergence when far from the optimum, if the objective function is strongly convex \cite{nedic2001convergence}. The rate analysis presented here fills this gap for a  class of convex problems that satisfy certain regularity conditions; see Sec. \ref{conv}.
	
	The use of dual subgradient algorithms for deterministic resource allocation was first popularized in \cite{kelly1997charging}. Recovery of near-optimum allocations via primal averaging was proposed in \cite{larsson1999ergodic}, and the result was extended to stochastic resource allocation problems in \cite{classconv,ale10}. The  corresponding convergence rate analysis for primal recovery was provided in \cite{nedic2009approximate}, which also serves as a starting point for the analysis presented here. Note however that the extension of the rate results to stochastic problems is not trivial, and does not follow immediately from the result in \cite{nedic2009approximate}. The specific assumptions required to develop the bounds in this paper are inspired from those used in the context of stochastic approximation and averaging \cite{Solo_1994}.
	
	From a broader perspective, the work in this paper is also related to the backpressure algorithm, first proposed in the context of stochastic network optimization \cite{neely2010stochastic}. As shown in \cite{yu2015convergence}, the dual subgradient algorithm when applied to deterministic resource allocation problems, may also be viewed as the so-called drift-plus-penalty algorithm. The analysis in \cite{yu2015convergence} however does not translate to convergence rate results for the {SDSD} algorithm.
	
	The wireless caching framework utilizing D2D communications was first proposed in \cite{ji2016wireless,archi_D2D}, where the fundamental limits were analyzed. The system model described here builds upon the basic framework of \cite{archi_D2D} by formulating the problem within the resource allocation fabric, and adding some implementation details. The results presented here may also be applied to other stochastic resource allocation formulations, such as those in broadcast power allocation\cite{stoch_res_prob_3},  OFDM \cite{ofdm}, beamforming \cite{sidiropoulos2006transmit}, cognitive radio networks \cite{marques3}, network utility maximization \cite{resource_alloca_tut_4,crolayketan,ale10}, demand-response in the smart grid \cite{gatsis2014stochastic}, \colb{smart grid powerded green communications \cite{wang2016dynamic,wang2016dynamic2}} and energy harvesting \colb{\cite{wang2016stochastic}}.
\vspace{-5mm}

\section{Problem Formulation}\label{probfor}
\colb{This section formulates the general stochastic network resource allocation problem. We begin with detailing a D2D caching example that is used to motivate the general problem.	}
\vspace{-6mm}	
	\subsection{\colb{Motivation: D2D Mobile Caching}}
	 The D2D framework enables direct communication between nearby user equipments (UE), enabling greater spectrum utilization, higher energy efficiency, and increased overall throughput. The technology also allows unique solutions to connectivity problems that arise at the network edge or as a result of cellular congestion at overcrowded events. As an example,  the D2D architecture proposed in \cite{archi_D2D} considers caching of popular content on mobile devices with excess storage. The content files are then available for download over a D2D link, allowing users to reach higher data rates, avoid congestion, and overcome coverage issues. By directly involving the smartphone equipped users into the process of content distribution, such an \emph{edge-caching} solution not only cuts down the hardware provisioning costs but also promises better user experience. 
	 
\begin{figure}
	 	\centering
	 	\includegraphics[scale=0.3]{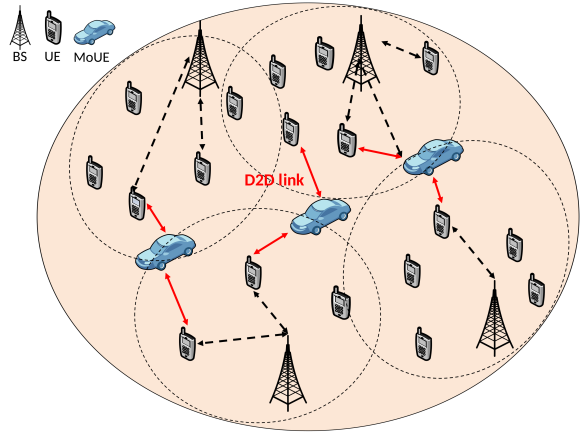}
	 	\caption{System model for mobile caching in D2D networks.}
	 	\label{fig0}
	 	\vspace{-1cm}
	 \end{figure}

	 This example builds upon the mobile caching framework proposed in \cite{archi_D2D}. Specifically, the mobile user equipment {(MoUE)} seeks to download a large file or stream media for a sufficiently long duration, while maintaining a reasonable download rate or quality of experience. {Let $\mathcal{M}\!\! = \!\{\!1, \ldots, M\!\}$ be the set of mobile caches in the network and at time $t$, the requested chunk be available at $\mathcal{M}_t \!\!\subset\! \!\mathcal{M}$ unique mobile caches (devices) that are at close proximity to the user.} The potential download rate $R_i(p_t^i, \!\gamma_t^i)$ depends on the power allocation $p_t^i$ at the $i$-th mobile cache, as well as on the channel gain $\gamma_t^i$ of the D2D links. The downloads also incur a cost $c_t^i\!\! \in \!\!\Rn_{++}$ per unit of transmit power $p_t^i$ for slot $t$. The costs could be in form of incentives provided to the mobile caches by the content delivery network (CDN) company, and/or directly charged to the MoUE in form of an ``enhanced coverage'' fee. At each time $t$, the {MoUE} selects a cache $i_t$ to download from, and obtains an average throughput of $r$ over time. Finally, the user satisfaction for the achieved average throughput $\!r$ is quantified through the concave utility function $U\!(r)$. Fig. \ref{fig0} depicts an example scenario, where a {MoUE} connects to different UEs in order to download cached data, that would otherwise be available only from the base stations. 
	 
	 The resulting stochastic resource allocation problem is formulated as
	 \vspace{-0mm}
	 \begin{subequations} \label{d2d0}
	 	\begin{align}
	 	\max_{r,\{p^i\}}~~ &U(r) - \Ex{{\sum_{{i\in\mathcal{M}_t}} c_{t}^ip_{t}^i}} \label{obj}\\
	 	\text{s. t. } & r\leq \Ex{{\sum_{{i\in\mathcal{M}_t}}R_i(p_{t}^i,\gamma_{t}^i)}} \label{ravg}\\
	 	&\{p^i\}_{{i\in\mathcal{M}}} \in \Pc, \ \ r_{\min} \leq r \leq r_{\max}
	 	\end{align}
	 \end{subequations}
	 where the expectations are with respect to the random variables $\bg_t\! :=\!{(\mathcal{M}_t, \{c^i_t\}_{i\in\mathcal{M}_t}, \{\gamma^i_t\}_{i\in\mathcal{M}_t})}$. The optimization variables consist of the power allocation function $p^i$ and the rate variable $r$. The formulation of  \eqref{d2d0} follows the classical ``utility minus penalty'' maximization format common to network resource allocation problems \cite{resource_alloca_tut_4}. The set of functions $\Pc$ is such that only one out of ${\{p_t^i\}_{i\in\mathcal{M}}}$ is non-zero for each $t$ (cf. Sec. \ref{d2d}). Consequently, the summations in \eqref{obj} and \eqref{ravg}  consist of only one term for each $t$. The set $\Pc$ also specifies the maximum and minimum values of $p_t^i$ for each $t$. Finally, the constraint in \eqref{ravg} ensures that the power allocated per-time slot is sufficient to satisfy the average rate requirement.
	 
	 The specific form of the rate function depends on the wireless technology used by the users. For instance, under slow fading scenarios, the power allocation and user selection can occur every coherence interval. Since channel state information can be acquired easily, the users may employ adaptive modulation and coding in order to achieve a rate close to the ergodic capacity. Specifically, for the mobile device $i$, the potential transmission rate is of the form \colb{$R_i(p_t^i,\gamma_t^i) : = W\colb{\log_2} (1 + p_t^i\gamma_t^i/\alpha)$ where $W$ is the bandwidth of the channel and $\alpha$ includes the effect of noise and interference as well as other impairments, such as the use of finite block length codes and imperfect channel state information at the transmitter.}
	 
	 More realistically, under fast fading scenarios, the power allocation must occur over intervals that are significantly longer than the coherence time. In this case, it is more reasonable to consider the average rate over several coherence intervals as $\!R_i(p_t^i,\!\gamma_t^i)\!:=\! W\Eh{\colb{\log_2}(1\!+\!p_t^i\gamma_t^ih_i/\alpha)}$ where $h_i$ is the small-scale fading gain \cite{tse}. In this case, the power allocation and user selection occur only on the basis of the average channel gain $\gamma_t^i$, which changes slowly. It is remarked that the system model considered here allows other forms of the rate function as well. 
	 \vspace{-7mm}
	\subsection{\colb{Stochastic Resource Allocation}}
	This section considers the more general stochastic resource allocation problem where the formulation involves expectations with respect to a collection of $q$ random variables with unknown distributions, {denoted by} $\bg \in \Rn^q$. Of particular interest are the problems arising in the context of wireless communications and networks, where $\bg$ captures the state of the system, and the  formulation takes the form \cite{crolayketan, ale10,marques3,ofdm} \vspace{-4mm}
\begin{align}
\hspace{-1cm}\pa \hspace{1cm}  (\x^\star, \p^\star) &= \argmax f_0(\x) \label{primal}\\
\text{s. t. }  & \u(\x) + \Ex{\v(\bg,\p_\bg)} \geq 0, \label{pc1}\\
& \x \in \mathcal{X}, \p \in \Pc. \label{pc2}
\end{align}
The optimization variables in $\pa$ comprise of the resource allocation variable $\x \in \Rn^d$ and the policy functions $\p:\Rn^q \rightarrow \Rn^p$. The objective function $f_0: \Rn^d \rightarrow \Rn$ is a concave function, while the set  $\mathcal{X}$ is {compact and }convex. The vector-valued constraint function is defined as $\u(\x):=[u_1(\x) \cdots u_K(\x)]^T$, where $\{u_i(\x):\Rn^d \rightarrow \Rn\}_{i=1}^K$ are concave functions.  In contrast, no such restriction is placed on the vector-valued function $\v:\Rn^q \times \Rn^p \rightarrow \Rn^K$ and the {compact} set of functions $\Pc$. The rate analysis in Sec. \ref{conv} however does require the overall problem to satisfy certain regularity properties, such as Slater's constraint qualification and differentiability of the subgradient error; see (\textbf{A1})-(\textbf{A4}). 

\colb{It can be seen that the D2D edge caching problem in \eqref{d2d0} is a special case of \eqref{primal}-\eqref{pc2}. Introducing a scalar variable $z \in \Rn$, it is possible to write \eqref{d2d0} equivalently as
	 \begin{subequations} \label{d2d011}
	 	\begin{align}
	 	\max_{r,\{p^i\}}~~ &U(r) +z  \label{obj11}\vspace{-5mm}\\
	 	\text{s. t. } &  -z-\Ex{{\sum_{{i\in\mathcal{M}_t}} c_{t}^ip_{t}^i}}\geq 0\\
	 	& -r+\Ex{{\sum_{{i\in\mathcal{M}_t}}R_i(p_{t}^i,\gamma_{t}^i)}} \geq 0\label{ravg11}\\
	 	&\{p^i\}_{{i\in\mathcal{M}}} \in \Pc, \ \ r_{\min} \leq r \leq r_{\max}
	 	\end{align}
	 \end{subequations}  
Comparing \eqref{d2d011} with \eqref{primal}-\eqref{pc2}, we see that $\x = [r, z]^T$ and $\bg_t :={(\mathcal{M}_t, \{c^i_t\}_{i\in\mathcal{M}_t}, \{\gamma^i_t\}_{i\in\mathcal{M}_t})}$. Likewise the forms of vector functions $\u$ and $\v$ can be readily inferred.} 
	
Since the distribution of $\bg$ is also not known in advance, it is generally not possible to solve $\pa$ in an offline manner. The goal here is to solve $\pa$ in an online fashion by observing the realizations of the independent identically distributed (i.i.d.) process $\{\bg_t\}_{t\in \Nn_0}$, where $\Nn_0$ is the set of non-negative integers. For most problems of interest, such a framework also entails online allocation of resources for each time $t$. To this end, the class of algorithms considered here will output the sequence of vector pairs $\{\x_t, \mathring{\p}_t\}$ for each $t$, for the purpose of allocating resources. For the sake of brevity, we will subsequently denote policy function $\p_t:= \p_{\bg_t}$ and $\s_t(\x,\p_t):=\u(\x) + \v(\bg_t,\p_{\bg_t})$, so that \eqref{pc1} can equivalently be written as $\Ex{\s_t(\x,\p_t)} \geq 0$. Here, it is understood that the expectation is with respect to the random vector $\bg_t$. Having introduced the problem at hand, we detail an example formulation in the context of D2D mobile caching.
	
%	\textbf{Mobile caching in D2D network:}
	\vspace{-5mm}
	\section{Solution via Dual Descent}\label{solu_dual}
	This section details the SDSD algorithm for solving $\pa$ in an online fashion. To this end, the basic assumptions are first stated (Sec. \ref{bass}), followed by the {SDSD} algorithm (Sec. \ref{sdda}), and a discussion of the known results (Sec. \ref{known}). 
	\vspace{-7mm}
	\subsection{Basic assumptions}\label{bass}
	The following assumptions are commonly utilized by different dual algorithms proposed in the literature. None of these assumptions are too restrictive, \colb{and they can be} easily verified  for most resource allocation problems of interest. 
	
	\begin{enumerate}
		\item[\textbf{\textbf{A1}}. ] \textbf{Non-atomic probability density function: } The random variable $\bg_t$ has a non-atomic probability density function (pdf). 
		\item[\textbf{A2}. ] \textbf{Slater's condition: } There exists strictly feasible ${(\tilde{\x},\tilde{\p})}$, i.e., ${\Ex{\st(\tilde{\x},\tilde{\p}_t)}} {>} 0$. 
		\item[\textbf{A3}. ] \textbf{Bounded subgradients: } The function $\st(\cdot,\cdot)$ takes bounded values, i.e., there exists a constant $G < \infty$ such that  $\norm{\s_t(\cdot,\cdot)}\leq G$ for all ${t\in\mathbb{N}_0}$. 
	\end{enumerate}
	
	In (\textbf{A1}), for $\bg_t$ to have a non-atomic pdf, it should not have any point masses or delta functions. Note that this requirement is not restrictive for a number of applications arising in wireless communications; see e.g. \cite{ale10}. The Slater's condition is also not restrictive, since a strictly feasible resource allocation can often be found for most real-world problems; see Sec. \ref{d2d}, \cite{crolayketan, ale10} for examples. Finally, the bound in (\textbf{A3}) also holds for most resource allocation problems where  the functions $\st(\cdot,\cdot)$ represent natural quantities such as instantaneous rate (cf. \eqref{ravg}), indicator function for channel outage \cite{marques3}, or household power consumption \cite{Gatsis_2}. Having introduced the basic assumptions, we are ready to state the {SDSD} algorithm. 
	\vspace{-6mm}
	\subsection{The {stochastic dual }subgradient algorithm}\label{sdda}
	Towards solving $\pa$, consider the more tractable dual formulation, which has a finite number of optimization variables. Introducing a dual variable $\lam \in \Rn^K_+$ corresponding to the constraint \eqref{pc1}, the Lagrangian is given by
	\vspace{-4mm}
	\begin{align}
	{L}(\lam,\x,\p)= f_0(\x) + \ip{{\lam,\Ex{{\st(\x,{\p}_t)}}}} \label{lag}
	\end{align} 
	where the constraints in \eqref{pc2} are kept implicit. The dual function is obtained by maximizing ${L}({\lam,\x,\p})$ subject to \eqref{pc2}, that is, 
		\vspace{-4mm}
	\begin{align}
	g(\lam) = \max_{\x \in \Xc, \p \in \Pc} {L}(\lam,\x,\p).
	\end{align}
	Finally, the dual problem of $\pa$ is given by
		\vspace{-4mm}
	\begin{align}
	\mathsf{D} = \min_{\lam \geq 0} &~ g(\lam).\label{dual} 
	\end{align}
	In general, since $\pa$ may be non-convex, it holds that $\mathsf{D} \geq \mathsf{P}$. It was shown in {\cite[Prop. 6]{crolayketan}}
	%\footnote{\colr{Proposition 6 establishes the zero duality gap for the stochastic constrained optimization problem \eqref{primal}-\eqref{pc2} considered in this paper.}} 
	however, that under (\textbf{A1})-(\textbf{A3}), it holds that $\mathsf{P} = \mathsf{D}$. The proof utilizes the Lyapunov convexity theorem, and holds even if at least one element of $\bg_t$ has an absolutely continuous cumulative distribution function (cdf) \cite{rao1972remark}. It is remarked that Lyapunov convexity has previously yielded similar results in control theory \cite{hermes1969functional}, economics \cite{shitovitz1973oligopoly}, and wireless communications \cite{luo2008dynamic,separation}. 
	
	The result on zero duality gap legitimizes the dual descent approach, since the dual problem is always convex, and the resultant dual solution can be used for primal recovery. To this end, similar problems in various contexts have been solved via the classical dual descent algorithm \cite{crolayketan, separation, classconv, marques3}, wherein the primal updates utilize various sampling techniques.
	
This paper considers the ergodic stochastic optimization (ESO) algorithm proposed in \cite{ale10} for a similar \colb{problem\footnote{\colb{The ESO algorithm is a stochastic dual subgradient descent algorithm applied to a resource allocation problem in \cite{ale10}.}}}. Applied to $\pa$, the ESO algorithm starts with an arbitrary $\lam_0$, and utilizes the following iterations for $t \in \mathbb{N}_0$,
\vspace{-4mm}
\begin{subequations}\label{eso}
	\begin{align} 
	(\x_t(\lamt),\p_t(\lamt)) &= \argmax_{\x \in \Xc, \mathring{\p} \in \Pi_t} f_0(\x) + {\ip{\lamt,\st({\x,\mathring{\p}})}}	\label{esop}\\
	\lam_{t+1} &= [\lamt - \epsilon \st(\x_t(\lamt),\p_t(\lamt))]_+. \label{esod}
	\end{align}
\end{subequations}
Here, $\Pi_t :=\{\p_{\bg_t} \in \Rn^p | \p \in \Pc\}$ is the set of all legitimate values of the vector $\p_{\bg_t}$ {and $\mathring{\p}$ denotes all feasible vectors $\p_t(\lamt)$}. The ESO algorithm is motivated from the fact that for any $\lam \in \Rn^K_{+}$, $\st(\x_t(\lam),\p_t(\lam))$ is a stochastic subgradient of the dual function $g(\lam)$. Consequently, the updates in \eqref{eso} amount to solving \eqref{dual} via the {SSD} algorithm with a constant step-size. The use of a constant step size is motivated from classical short memory adaptive algorithms such as the least mean squares algorithm. As stated earlier, the constant step-size algorithms can even handle abrupt changes in the problem parameters, without being restarted.
	
	This paper considers the projected variant of the {SSD} algorithm for dual updates. Specifically, the updates in \eqref{esod} are projected on to a compact set $\mathcal{L} \subset \Rn_{+}^K$, and take the following form
	\begin{align}
	\lam_{t+1} &= {\mathcal{P}_{\mathcal{L}}\left({\lamt - \epsilon \st(\x_t(\lamt),\p_t(\lamt))}\right) }\label{dualupp}
	\end{align}
	where for any $\lam \in \Rn^K$, 
	\vspace{-5mm}  
	\begin{align}
	[{\mathcal{P}_{\mathcal{L}}}\left({\lam}\right)]_i = \begin{cases} 0 & [\lam]_i < 0 \\
	\lambda_{\max} & [\lam]_i > \lambda_{\max} \\
	\lambda_i & 0 \leq [\lam]_i \leq \lambda_{\max}
	\end{cases}
	\end{align}
	for $1\leq i\leq K$. In other words, large values of $[{\lam}]_i$ are truncated to ${\lambda_{\max}}$, where $\lambda_{\max} \gg \norm{\ls}_{\infty}$. Such a modification is already applicable to any practical implementation of the SSD algorithm, where $\lam_t$ is not allowed to take arbitrarily large values. Since $\ls$ is not known in advance, a bound on $\norm{\ls}_{\infty}$ is derived in Appendix \ref{lambda} using (\textbf{A2}). Consequently, the following rule can be used for choosing $\lambda_{\max}$ in practice: 
	\vspace{-4mm}
	{\begin{align}\label{lambound}
		%\lambda_{\max} >> \frac{1}{\gamma} (g-f) >= \frac{g-f}{G}\\
		\lambda_{\max} \gg \frac{1}{{{\chi(\tilde{\x},\tilde{\p})}}}(g(\tilde{\lam}) - f_0(\tilde{\x})) \geq \frac{(g(\tilde{\lam}) - f_0(\tilde{\x}))}{G}
		\end{align} 
		where $\tilde{\lam} \in \Rn^K_{+}$, $({\tilde{\x}, \tilde{\p}})$ is a strictly feasible solution to $\pa$ (cf. (\textbf{A2})), {$G$ is the subgradient bound (cf. \textbf{A3})} and ${\chi(\tilde{\x},\tilde{\p})} := \min_{1\leq k\leq K} \Ex{[\st({\tilde{\x},\tilde{\p}_t})]_k}$.} {In general, the quantity $\chi(\tilde{\x},\tilde{\p})$ may be calculated empirically. However, for many problems of interest, a bound on $\lambda_{\max}$ may arise naturally (cf. Sec. VI).} \colb{The projected SSD proposed in \eqref{dualupp} ensures that the iterates $\lam_t$ stay bounded for all ${t \in \mathbb{N}_0}$. The boundedness condition is required for carrying out the rate analysis in Sec. \ref{conv}.}
		
	\vspace{-5mm}
	\subsection{Known results}\label{known}
	The asymptotic properties of the {SDSD} algorithm with constant step-size are well-known \cite{nedic2001convergence,bertsekas2011incremental}. Asymptotic convergence results for the ESO algorithm, applied to slightly different resource allocation problem, were established in \cite{ale10}. The results in \cite{ale10} can readily be extended to $\pa$ solved via projected {SDSD}, and take following form:
	\begin{subequations}\label{ceso}
		\begin{align}
		&\lim\limits_{{{t}} \rightarrow \infty} \frac{1}{{{t}}} \sum_{{{\tau}}=0}^{{{t}}-1}  \s_{{\tau}}(\x_\tau(\lam_\tau),\p_\tau(\lam_\tau)) \geq 0	& \text{a. s.}\label{ceso1}\\
		&\lim\limits_{{{t}} \rightarrow \infty} f(\bar{\x}_{{t}}) \ge \mathsf{P}-({\epsilon \bar{G^2}}/{2}) & \text{a. s.}\label{ceso2}
		\end{align}
	\end{subequations}
	where the running average $\bar{\x}_t$ is defined as 
	\vspace{-2mm}
%	\begin{align}
$	\bar{\x}_{{t}} = \frac{1}{{t}} \sum_{{\tau}=0}^{{t}-1}  \x_\tau(\lam_\tau) \label{run_avg_x}$
%	\vspace{-5mm}
%	\end{align}
	and $\bar{G^2}$ is the bound on $\Ex{\norm{\st(\x_t(\lam_t), \!\p_t(\!\lam_t\!))}^2\!}$. An important feature of the stochastic algorithm is that the primal updates in \eqref{esop} can be used for allocating resources in real-time. Further, such allocations will be asymptotically feasible and near-optimal for almost every realization of the random process $\{\bg_t\}_{{t\in \Nn_0}}$. 
	\vspace{-10mm}
\section{Convergence Rate Results}\label{conv}
This section develops various results regarding the rate of convergence of the SDSD algorithm. In contrast to the asymptotic results in \eqref{ceso}, the goal here is to quantify the rate at which the allocations specified by \eqref{esop} become optimal. Such results are of practical significance to the protocol designers, since they can be used to estimate the number of iterations required for the primal and dual objectives to be near-optimal. In the case of the constant step-size {SDSD}, the convergence rate also depends on the step-size parameter $\epsilon$. For instance, it is well-known that the choice $\epsilon \rightarrow 0$, motivated from the result in \eqref{ceso2}, leads to slow convergence in all constant step-size (sub-)gradient descent algorithms. The results presented here provide a precise characterization of the trade-off between $\epsilon$ and the convergence rate for the {SDSD} algorithm.
	
As in \cite{ale10}, the results in this section make use of the strong law of large numbers, and thus hold for almost every realization of the i.i.d. process $\{\bg_t\}_{{t\in \Nn_0}}$. It is emphasized that the analysis presented here is quite different from the standard convergence analysis carried out for SSD algorithm and its variants \cite{nedic2001convergence,bertsekas2011incremental,bottou_91c_neural_learning,bottou2010large}. It is also different from the {non-asymptotic} analysis for the case of diminishing step-size SSD algorithms, that only applies to ensemble averages \cite{bach2011non}. Furthermore, the rate results presented in \cite{bach2011non} apply only to the unconstrained stochastic subgradient algorithm, and cannot be extended to constrained problems (cf. $\pa$) or to the projected subgradient algorithm (cf. \eqref{dualupp}).
	
The results are first developed for the general SSD algorithm (Sec. \ref{ssdrate}), and subsequently specialized to the resource allocation problem at hand (Sec. \ref{sddrate}). 
\vspace{-4mm}
\subsection{Convergence rate for the SSD algorithm}\label{ssdrate}
This section considers the generic optimization problem 
\vspace{-4mm}
\begin{align}
\lam^\star &= \arg\min_{\lam \in \Lambda} g(\lam) \label{generic}
\end{align}
where $\Lambda \subset \Rn^K$ is a closed, compact, and convex set, and $\max_{\lam \in \Lambda}\norm{\lam-\lam^\star} \leq \Lambda_{\max} < \infty$. Similar to \eqref{dual}, the optimum function value is denoted by $\mathsf{D} = g(\lam^\star)$. Given $\lam \in \Lambda$, let  $\f(\lam) \in \partial g(\lam)$ be a subgradient of $g(\lam)$. Similar to $\pa$, let $\ft(\lam):=\f_{\bg_t}(\lam)$ for all ${t\in \Nn_0}$ be the corresponding stochastic subgradients that depend on the i.i.d. process $\{\bg_t\}_{{t\in \Nn_0}}$ and satisfy $\Ex{\ft(\lam)} = \f(\lam)$ for any $\lam \in \Lambda$. For instance, in the simplest case, the stochastic gradient could be of the form $\ft(\lam)=\f(\lam)+\boldsymbol{\zeta}_t$, where $\boldsymbol{\zeta}_t$ is a zero mean i.i.d. random variable. The optimization problem \eqref{generic} is solved via the projected SSD algorithm. 
\begin{align}\label{ssdup}
\lam_{t+1} &= \pl{\lamt - \epsilon \ft(\lamt)} 
\end{align}
where $\pl{\cdot}$ denotes the projection operation. The algorithm is initialized with an arbitrary $\lam_0 \in \Lambda$ such that $B_0:=\norm{\lam_0-\lam^\star} < \infty$. Next, we make certain assumptions specific to \eqref{generic}. To this end, define the stochastic error $\e_t(\lam):= \f_t(\lam)-\f(\lam)$, and observe that for any ${\lam \in \Lambda}$, the sequence $\{\e_t(\lam)\}$ is also i.i.d.
	
	\begin{enumerate} 
		\item[\textbf{\textbf{A3}}$^\prime\hspace{-1mm}.$ ] \textbf{Bounded subgradients: } There exists constant $G < \infty$ such that $\norm{\f_t(\lam)} \leq G$ for all $\lam \in \Lambda$. 
		\item[\textbf{A4}. ] \textbf{Continuously differentiable error: } The error function $\e_t(\lam)$ is continuously differentiable on $\Lambda$, and the gradient with respect to $\lam$ satisfies $\norm{\nabla_{\lam}\e_t(\lam)} < G_e < \infty$. 
	\end{enumerate}
		
	The requirement for bounded subgradient in (\textbf{A3}$^\prime$) is analogous to that in (\textbf{A3}) for $\pa$. \colb{Here, (\textbf{A3}$^\prime$) is stated separately because the problem in \eqref{generic} is more general than the dual of \eqref{primal}-\eqref{pc2}. In practice, applying the results of this section to \eqref{dual} entails substituting $\f_t(\lambda)=\s_t(\x(\lam),\p_t(\lam)$, which makes (\textbf{A3'}) the same as (\textbf{A3}).} The error function $\e_t(\lam)$ may not always be continuous or differentiable for the problem at hand, and the same must be verified explicitly. It is emphasized that (\textbf{A4}) need only be checked for $\e_t(\lam)$ and not for $\f_t(\lam)$, which is still allowed to be non-differentiable; see Sec. \ref{d2d}.
		
		As an example, consider the class of problems where the convex objective function takes the form $g(\lam) \!=\! \Ex{\ell_t(\lam)}\! +\! r(\lam)$, where $\ell_t(\!\lam\!)$ is a twice-differentiable loss function that depends on the `data index' $t$, and $r(\lam)$ is a possibly non-differentiable regularizer. For such problems, the error function becomes $\e_t(\lam) \!\!= \!\!\nabla_{\lam} \ell_t(\lam) \!-\! \Ex{\nabla_{\lam} \ell_t(\lam)}$, which is clearly differentiable. Further, the `loss-plus-regularizer' problem structure is quite general, and includes well-known formulations such as LASSO \cite{bottou2010large} and nuclear norm regularized matrix least squares \cite{toh2010accelerated}. Specifically, given regressands $\{\!y_t\!\}$ and regressors $\{\!\x_t\!\}$, the objective function in the LASSO formulation takes the form $\sum_t (y_t \!-\! \ip{\x_t,\!\lam}) \!\!+\!\! \norm{\lam}_1$ and thus adheres to (\textbf{A4}). The first result is regarding the objective function values obtained from \eqref{ssdup}, and holds for all $T$.

		%{For brevity, a simple example is elaborated for which the assumption in (\textbf{A4}) will always be satisfied. Considering the structure of constraint function presented in $\pa$, we can write the stochastic subgradient of the dual function $g(\lam)$ as}
		%{\begin{align}
		%\f_t(\lam)=\u(\x_{t}(\lam)) + {\v(\bg_t,\p_{\bg_{t}}((\lam)))}
		%\end{align}}
		%{Note that the error function $\e_t(\lam):=\f_t(\lam)-\Ex{\f_t(\lam)}=\v(\bg_t,\p_{\bg_{t}})-\Ex{\v(\bg_t,\p_{\bg_{t}})}$ will always satisfies (\textbf{A4}) if the function $\v(\cdot)$ turns-out to be linear function of $\lam$, because in that case $\e_t(\lam)$ will not depend upon $\lam$ and gradient will identically be equal to $0$. }
		%{\begin{align}
		%\e_t(\w)=
		%\begin{cases}
		%\w_t-\Ex{\w_t}, \text{if} \ \ y_t\ip{\w_t,\x_t} >1\\
		%\w_t-\lambda(y_t\x_t)-\Ex{\w_t-\lambda(y_t\x_t)}
		%\end{cases}
		%\end{align}which is essentially continuosly differentiable with respect to $\w_t$}.

\begin{thm} \label{conrate}
	Under (\textbf{A1})-(\textbf{A4}) and for $T = n/\epsilon$, the minimum dual function value is bounded as 
	\begin{align}\label{conrateeq}
	\min_{0\leq \tau \leq T-1} g (\lam_\tau) \leq \frac{1}{T} \sum_{\tau = 0}^{T-1} g(\lam_\tau) &\leq \mathsf{D} + \frac{B_0}{2n} + \frac{\epsilon G^2}{2} + C_T(n,\epsilon) 
	\end{align}
	where the random variable $C_T(n,\epsilon)$ \colb{holds} for $\zeta > 0$, 
	\begin{subequations}\label{ctlim}
		\begin{align}
		&\epsilon^{\zeta-1/2} C_T(n,\epsilon) \rightarrow 0 \hspace{.5cm}\text{a.s. as }\epsilon \rightarrow 0 \text{ for fixed }n < \infty \label{ctlimep}\\
		&n^{1/2-\zeta} C_T(n,\epsilon) \rightarrow 0  \hspace{.5cm}\text{a.s. as } n \rightarrow \infty \text{ for fixed }\epsilon > 0 \label{ctlimn}. \\
		&\mathbb{P}(C_t(n,\epsilon) > \nu^{\zeta-1/2}) < A\exp\left(-\nu^{2\zeta}\right) \label{ctprob}
		\end{align}
	\end{subequations}
	where $\nu :=\max\{\frac{1}{\epsilon},n\}$ and $A < \infty$ is a constant that does not depend on $n$ or $\epsilon$. 
\end{thm}

%\begin{thm} \label{conrate}
%\colr{Under (\textbf{A1})-(\textbf{A4}), the minimum dual function value is bounded as 
%\begin{align}\label{conrateeq}
%\hspace{-5mm}\min_{0\leq \tau \leq T-1} g (\lam_\tau) \leq \frac{1}{T} \!\sum_{\tau = 0}^{T-1} \!g(\lam_\tau) &\!\leq \!\mathsf{D}\! + \frac{G\norm{\lam_0}}{\sqrt{T}}  + C_T 
%\end{align}
%where the random variable $C_T$ is such that for any $\zeta > 0$, 
%\begin{subequations}\label{ctlim}
%\begin{align}
%T^{\frac{1}{4}-\zeta} C_T &\rightarrow 0 & \text{a.s. as } T&\rightarrow \infty 
%\end{align}
%\end{subequations}
%Further, given arbitrary $\nu > 0$, there exists a constant $A_\nu \in (0,1)$ not dependent on $T$ such that 
%\begin{align}\label{ctprob}
%\mathbb{P}(C_T > \frac{\nu}{\sqrt{T}^{1-\zeta}}) < A_\nu \exp(-T^\zeta)  
%\end{align}}
%\end{thm}}

		It is remarked that since $g(\cdot)$ is convex, the bound in \eqref{conrateeq} also holds for $g(\bar{\lam}_t)$, where $\bar{\lam}_t : = \frac{1}{t}\sum_{\tau = 0}^{t-1} \lam_\tau$ is the running average of the iterates. Theorem 1 characterizes the manner in which the minimum objective function value approaches $\mathsf{D}$ for large $t$. Of the three terms in this optimality gap, the first one depends on the initialization and decays as $\O(1/n)$. The second term depends on the subgradient bound, and decays linearly with the step-size $\epsilon$. Finally, the third term is random, and decays almost surely as $\O((\epsilon/n)^{1/2-\zeta})$ for any $\zeta >0$ (cf. \eqref{ctlim}). Alternatively, the probability of the third term being non-zero decays exponentially as either $n \rightarrow \infty$ or $\epsilon \rightarrow 0$ (cf. \eqref{ctprob}). Indeed, for a given run of \eqref{ssdup} with a fixed $\epsilon$, the probability of the third term being non-negligible starts to decrease only beyond $n > 1/\epsilon$ or equivalently, $T > 1/\epsilon^2$. 
			
			Further intuition on the convergence rate can be obtained by considering the two cases in  \eqref{ctlim}. When $\epsilon > 0$ is fixed, it can be seen that the asymptotic results in \cite{nedic2001convergence,bertsekas2011incremental, ale10} follow directly from Theorem \ref{conrate} as $n \rightarrow \infty$. That is, while the initial condition is ``forgotten'' for $t \gg 1/\epsilon$, the optimality gap does not necessarily approach zero, but is eventually bounded by $\epsilon G^2/2$. At the same time, the fluctuations due to the stochastic term subside exponentially fast; see \eqref{ctprob}. 
			
			On the other extreme, consider the case when $n$ is kept fixed, while the algorithm is run for different values of $\epsilon$. For the scenarios when $\epsilon$ is arbitrarily small, the asymptotic optimality gap is clearly negligible. However, for such small step-sizes, the algorithm takes a long time to forget the initial conditions, since the first term decays only as $\O(\frac{1}{\epsilon t})$. Consequently, for all runs when $\epsilon$ is taken to be small, the algorithm will appear to converge slowly. Likewise, the probability of the stochastic term being non-negligible starts to decrease exponentially only for $T > 1/\epsilon^2$ (cf. \eqref{ctprob}). It is remarked that such a trade-off also applies to the classical subgradient method \cite{nedic2009approximate}, and the result in Theorem 1 can be viewed as its stochastic counterpart. It is remarked that the results in \cite{nedic2009approximate} can be readily obtained by taking expectation on both sides of \eqref{conrateeq} since we have that $\Ex{C_t(n,\epsilon)}\!=\!0$. Observe further that unlike the results in \cite{bottou1998online,bottou2012stochastic,kushner1994analysis,bach2011non,recht2011hogwild} that hold on an average, the almost sure results in \eqref{conrateeq} cannot be specified in terms of problem parameters alone. Indeed, while it holds that $C_t(n,\eps) \leq 2G^2(\Lambda_{\max}+1)$, such a bound is not very useful in the present case, as compared to the stronger convergence rate result in \eqref{conrateeq}.
			
Finally, it is remarked that it may be possible to minimize the bound in \eqref{conrateeq} to the extent possible, by fixing $T$ and choosing a corresponding step size. In the present case, given $T$, the bound is the smallest when $\epsilon = 1/\sqrt{T}$ which yields the following result
\begin{align}
\min_{0\leq \tau \leq T-1} g (\lam_\tau) \leq \frac{1}{T} \!\sum_{\tau = 0}^{T-1} \!g(\lam_\tau) &\!\leq \!\mathsf{D}\! + \mathcal{O}\left(\frac{1}{\sqrt{T}}\right)  + C_T 
\end{align}
where the random variable $C_T = \mathcal{O}(T^{-1/4})$ almost surely. The result in Theorem 1 may therefore be seen as the generalization of the results in \cite{bottou2010large,bach2011non} that have also reported an $\mathcal{O}(T^{-1/2})$ bound on average but have not analyzed the almost sure behavior. It is emphasized however that in practice, minimizing the bound may not necessarily translate to an improved convergence rate. Moreover, the number of iterations $T$ for which the algorithm runs may not necessarily be known in advance, e.g., in target tracking applications. Instead, it may be simpler to specify a fixed value of $\epsilon$, and continue to run the algorithm till the contribution of the $\mathcal{O}(\frac{1}{n})$ term becomes tolerably small. 
			 
		Before proceeding with the proof of Theorem \ref{conrate}, an intermediate lemma establishing rate results on various time-averages is provided. The proof of Theorem \ref{conrate} will subsequently utilize these results by expressing the optimality gap in \eqref{conrateeq} in terms of these time-averages. 
			
			\begin{lem} \label{errbounds}
				Let $\T:=\{t_1, t_2, \ldots, t_{\abs{\T}}\}$ be a set of natural numbers such that $t_i \neq t_j$. Then for any $T \geq \abs{\T}$ and $\lam, \lam' \in \Lambda$, it holds under {(\textbf{A3}$^\prime$)-(\textbf{A4})} that
				\small
				\begin{subequations}\label{a5a6}
					\begin{align}
					\norm{\frac{1}{T} \sum_{t \in \T} \e_t(\lam)} &\leq L^1_{T}(\T)\\
					\norm{\frac{1}{T} \sum_{t \in \T} \e_t(\lam) - \e_t(\lam')} &\leq L^2_{T}(\T)\norm{\lam-\lam'} 
					\end{align}
				\end{subequations}
			\normalsize
				where, for a given $\zeta > 0$,  the random variables $\{{L^i_{T}(\T)}\}_{i=1,2}$ satisfy 
				\begin{subequations} \label{sllnrate}
					\begin{align}
					T^{1/2-\zeta}L^i_T(\T) &\rightarrow 0 & \hspace{-5mm}\textrm{a.s. as } T \rightarrow \infty \label{sllnrate1}\\
					\mathbb{P}(L^i_T(\T) > T^{\zeta-1/2}) &< A_i \exp(-T^{2\zeta}) \label{sllnrate2}
					\end{align}
				\end{subequations}
				where the constant $A_i < \infty$ does not depend on $T$. 
			\end{lem}
				
				\begin{IEEEproof} {Observe that the i.i.d. process $\{\e_t(\lam)\}_{t \in \T}$ is zero-mean and satisfies
						\vspace{-3mm}
						\begin{align}
						\norm{\e_t(\lam)} \leq \norm{\f_t(\lam)} + \norm{\f(\lam)} \leq 2G
						\end{align}
						where the last inequality holds from (\textbf{A4}). Therefore, it follows from the strong law of large numbers that for any $T \geq \abs{\T}$, 
					%	\begin{align}
					$	\frac{1}{T}\sum_{t\in \T} \e_t(\lam) \rightarrow \mathbf{0}$
					%	\end{align}
						almost surely as $T \rightarrow \infty$. It can also be seen that the same holds for $L_T^1(\T) := \norm{\frac{1}{T}\sum_{t\in \T} \e_t(\lam)}$. The rate results in \eqref{sllnrate} hold as consequences of the strong law of large numbers for i.i.d. sequences with bounded moments; see {\cite[Chap. 7]{slln}} for \eqref{sllnrate1}. {Finally, \eqref{sllnrate2}, follows from the Bernstein inequality  applied to i.i.d. zero-mean and bounded random vectors $\{e_t(\lam)\}$ \cite{fuk1971probability}.}}
					
				Denote the $j$-th entry of $\e_t(\lam)$ by $e^j_t(\!\lam\!)$ for $1 \!\leq\! j \!\leq K$. From {(\textbf{A3}$^\prime$)-(\textbf{A4})}, we have that $e^j_t(\lam)$ is bounded and continuously differentiable on $\Lambda$. Consider arbitrary $\lam \!\neq \!\lam'\! \in\! \Lambda$, and observe that since $\Lambda$ is convex, it holds for any $\beta\! \in\! [0,\!1]$ that $\lam_{\beta}\!:=\!\beta\lam \!+\! (\!1\!-\!\beta)\lam' \!\in \!\Lambda$. It is now possible to use the mean-value theorem, which guarantees that there exists some $\beta_j \!\in\! [0,\!1]$, such that
				\vspace{-4mm}
				\begin{align}\label{mvt}
				e^j_t(\lam) - e^j_t(\lam') = \langle \nabla e^j_t(\lam_{\beta_j}), \lam-\lam' \rangle.
				\end{align}
				Here, $\nabla \e^j_t(\lam_{\beta_j})$ is an i.i.d. random variable that is also zero-mean, since for continuously differentiable and bounded functions { (cf. (\textbf{A4}))}, we have that $\Ex{\nabla \e^j_t(\lam_{\beta_j})} = \nabla \Ex{e^j_t(\lam_{\beta_j})} = 0$. Taking summation in \eqref{mvt} and stacking the $K$ components, it follows for any $T \geq \abs{\T}$, that
				\begin{align}\label{matrix1}
				\frac{1}{T}\sum_{t\in \T} \e_t(\lam) - \e_t(\lam') = \E_T(\lam,\lam')(\lam-\lam')
				\end{align}
				where the $K \times K$ matrix $\E_T(\lam,\lam')$ is defined as $[\E_T(\lam,\lam')]_{jk} : = {\frac{1}{T}\sum\limits_{t\in \T}}[\nabla e^j_t(\lam_{\beta_j})]_k$, {where the subscript is used to denote the $k$-{th} element of vector $\nabla e^j_t(\lam_{\beta_j})$}. Applying the Cauchy-Schwarz inequality to \eqref{matrix1}, we obtain
			%	\begin{align}
				$\norm{\frac{1}{T}\sum_{t\in \T} \e_t(\lam) - \e_t(\lam')} \leq \norm{\E_T(\lam,\lam')}\norm{\lam-\lam'}$.
			%	\end{align} 
				From the strong law of large numbers, we have that $[\E_T(\lam,\lam')]_{jk} \rightarrow 0$ almost surely as $T \rightarrow \infty$ for all $1\leq j,k\leq K$. It can be seen that the same also holds for $L^2_T(\T):=\norm{\E_T(\lam,\lam')}$. Finally, the rate results in \eqref{sllnrate} {follow {from} \cite[Chap. 7]{slln}$^2$ and the Bernstein inequality  applied to i.i.d. zero-mean and bounded random variables $\{[\nabla e_t^j(\lam_{\beta_j})]_k\}$ \cite{fuk1971probability}. }
				\end{IEEEproof}
				
				The proof of Theorem \ref{conrate} follows in two steps: the derivation of the overall form required in \eqref{conrateeq}, presented next; and the analysis of the random term $C_t(n,\epsilon)$ deferred to Appendix \ref{appbound}.

			\begin{IEEEproof}[Proof of Theorem \ref{conrate}]
				In order to derive the bound in \eqref{conrateeq}, recall that since $g(\lam)$ is convex, we have that,
				\vspace{-8mm}
				\begin{align}
				g(\lamt) &\leq g(\ls) + {\ip{\fa(\lamt),\lamt - \ls}} & t&\in \mathbb{N}_0. \label{convexg}
				\end{align}
				Letting $g_t := g(\lamt)$, it follows from the non-expansive property of $\pl{\cdot}$ that
				\vspace{-3mm}
				\small
				\begin{align}
				\hspace{-3.2mm}\delta_{t+1} &:= \norm{\lam_{t+1}-\ls}^2 = \norm{\pl{\lamt - \mu \ft(\lamt)}-\ls}^2 \label{delta} \\
				&\leq \norm{\lamt-\ls - \mu\ft(\lamt)}^2 \label{nonexp}=  \norm{\lamt-\ls}^2 - 2\mu{\ip{\ft(\lamt),\lamt-\ls}} + \mu^2\norm{\ft(\lamt)}^2 \\
				&\leq \delta_t - 2\mu{\ip{\fa(\lamt),\lamt-\ls}} + \mu^2G^2- 2\mu{\ip{\ft(\lamt)- \fa(\lamt),\lamt-\ls}} \label{bsub}\\
				&\leq \delta_t - 2\mu \left(g_t-\mathsf{D}\right) -  2\mu{\ip{\ft(\lamt)-\fa(\lamt),\lamt-\ls}}+ \mu^2G^2 \label{conv1}
				\end{align}
				\normalsize
				where \eqref{bsub} follows from ({\textbf{A3'}}) and \eqref{conv1} follows from \eqref{convexg}. Rearranging  \eqref{conv1} yields
				\begin{align}
				2\mu \left(g_t-\mathsf{D}\right) \leq & \left(\delta_t -\delta_{t+1}\right)  -  2\mu{\ip{\ft(\lamt)-\fa(\lamt),  \lamt-\ls}}+ \mu^2G^2
				\end{align}
				Taking sum over $\tau = 0, 1, \ldots, t-1$ and noting that $B_0=\delta_0$ and that $\mu t \geq n$, yields 
				\begin{align}
				\frac{1}{t}\!\sum_{\tau = 0}^{t-1} \!{g_\tau} &\!\leq \!\mathsf{D} \!+\! \frac{\delta_0 \!-\! \delta_{t}}{2\mu t}\!-\! \frac{1}{ t} \!\sum_{\tau = 0}^{t-1} {\ip{\f_\tau(\lam_\tau)\!-\fa(\lam_\tau),\lam_\tau\!-\!\ls}}\!+\! \frac{\mu G^2}{2}\!\leq\! \mathsf{D}\! +\! \frac{B_0}{2n} \!+\! \frac{\mu G^2}{2} \!+\! C_t(n,\!\epsilon\!) \label{bound1}
				\end{align}
				where the last inequality follows since $\delta_{t} \geq 0$ and the stochastic term in \eqref{bound1} is defined as 
				\begin{align} \label{ctne}
				C_t(n,\epsilon) := \abs{\frac{1}{ t}\sum_{\tau = 0}^{t-1} {\ip{\f_\tau(\lam_\tau)-\fa(\lam_\tau),\lam_\tau-\ls}}}.
				\end{align}
				Since \eqref{bound1} is of the same form as required in \eqref{conrateeq}, it remains to show that $C_t(n,\epsilon)$ converges in the sense of \eqref{ctlim}-\eqref{ctprob}. The convergence analysis for $C_t(n,\epsilon)$ makes use of the bounds developed in Lemma \ref{errbounds} and is deferred to Appendix \ref{appbound}. 
			\end{IEEEproof}
					
				It is remarked that the results in Theorem \ref{conrate} can likely be generalized to the case when $\Lambda$ is not necessarily compact. Such a generalization is likely possible because the strong law of large numbers, as well as the rate results in {\cite[Chap. 7]{slln}} and {\cite[Theorem 1]{baum1962exponential}}
				%This result states that for an exponentially converging stationary sequence of random variables $\{X_t\}$, there exists a $t_0$ such that $\Ex{\exp(tX_t)}<\infty$ for $t\in						[-t_0,t_0]$.} } 
				only require the random process to have bounded moments. Nevertheless, the requirement that $\norm{\lam_t-\ls} \leq \Lambda_{\max} < \infty$ is not too restrictive, and greatly simplifies the analysis.
				\vspace{-5mm}
					\subsection{Convergence rate for the {SDSD} algorithm}\label{sddrate}
					In order to apply the results developed in Sec. \ref{ssdrate} to the dual problem \eqref{dual}, observe that the stochastic subgradient of $g(\lam)$ for any $\lam \in \Rn^{K}_{+}$ is given by 
					\begin{align}
					\ft(\lam) &= \st(\x_t(\lam),\p_t(\lam)) \label{ftlamdual}\\
					\x_t(\lam) &:= \argmax_{\x \in \Xc} f_0(\x) + {\ip{\lam,\u(\x)}}\\
					\p_t(\lam) &:= \argmax_{\mathring{\p} \in \Pi_t} {\ip{\lam,\v(\bg_t,\mathring{\p})}}. \label{pup}
					\end{align}
					With $\ft(\lam)$ as defined in \eqref{ftlamdual}, the projected SSD updates take the same form as \eqref{ssdup}, with $\Lambda_{\max} = 2\sqrt{K}\lambda_{\max}$. Further the bound required in {(\textbf{A3}$^\prime$)} follows from {(\textbf{A3})}. Therefore, Theorem \ref{conrate} applies as is to the dual objective function under (\textbf{A3}) and {(\textbf{A4})}. 
					
					For the resource allocation problem however, the behavior of the primal objective function is more important. The subsequent theorem characterizes the primal near-optimality when the running average of $\{\x_t(\lam_t)\}$ is used for allocating resources. For the purpose of rate analysis, time is divided into epochs of duration $\frac{1}{\epsilon}$ each, and the result is expressed in terms of $\epsilon$ and $n$.
				\vspace{-5mm}	
\begin{thm} \label{conv_primal}
	Under (\textbf{A1})-(\textbf{A4}),		and for $n/\epsilon \leq t < (n+1)/\epsilon$, the average primal objective function is near optimal in the following sense: 
	\begin{align}
	f_0(\bar{\x}_t) \geq \frac{1}{t}\sum_{\tau = 0}^{t-1}f_0(\x_t) \geq \mathsf{P} - \frac{R_0}{2n} - \frac{\epsilon G^2}{2} - {C'_t(n,\epsilon)} \label{t2}
	\end{align}
	where $R_0:=\norm{\lam_0}^2$, $\bar{\x}_t = \frac{1}{t}\sum_{\tau = 0}^{t-1} \x_\tau$, and the random variable ${C'_t(n,\epsilon)}$ is such that for $\zeta > 0$, 
	\begin{subequations}\label{ctlim2}
		\begin{align}
		&\epsilon^{\zeta-1/2} {C'_t(n,\epsilon)} \rightarrow 0 \ \ \  \text{a.s. as }\epsilon \rightarrow 0 \text{ for fixed }n < \infty \label{ctlimep2}\\
		&n^{1/2-\zeta} {C'_t(n,\epsilon)} \rightarrow 0  \ \ \  \text{a.s. as } n \rightarrow \infty \text{ for fixed }\epsilon > 0 \label{ctlimn2} \\
		&\mathbb{P}(C_t(n,\epsilon) > \nu^{\zeta-1/2}) < A\exp\left(-\nu^{2\zeta}\right)\label{ctprob2}
		\end{align}
	\end{subequations}
	where $\nu :=\max\{\frac{1}{\epsilon},n\}$ and $A < \infty$ is a constant that does not depend on $n$ or $\epsilon$.
\end{thm}
						
					The term  $C'_t(n,\epsilon)$ in Theorem \ref{conv_primal} is very similar to $C_t(n,\epsilon)$ in Theorem \ref{conrate}, and therefore decays at the same rate. It follows from Theorem \ref{conv_primal} that the resource allocation yielded by the projected {  SDSD} algorithm is near optimal since the average primal objective value is close to $\mathsf{P}$. Similar to \eqref{conrateeq}, the bound in \eqref{t2} also holds for ${\max_{0 \leq \tau \leq t-1} f_0(\x_\tau)}$, as well as for $f_0(\bar{\x}_t)$. Further, the optimality gap in \eqref{t2} is also similar to the one in \eqref{conrateeq}, and therefore decays at the same rate. For details, see the discussion after the statement of Theorem \ref{conrate}. 
						
					In order to prove Theorem \ref{conv_primal}, the specific form of the bound in \eqref{t2} is first established. The rest of the proof is much the same as before, and results from Lemma \ref{errbounds} are again used to derive the bounds on $C'_t(n,\epsilon)$ as in the proof of Theorem \ref{conrate}.

					\begin{IEEEproof}[Proof of Theorem \ref{conv_primal}]
						Recall that the subgradient of $g(\lam)$ is given by $\f(\lam) = \Ex{\st(\x_t(\lam), \p_t(\lam))}$, so that $g(\lam)=f_0(\x(\lam)) + {\ip{\lam,\f(\lam)}}$. Since $f_0$ is concave, the following inequalities hold:
						\small
						\begin{align}
						f_0(\bar{\x}_t) &\geq \frac{1}{t}\sum_{\tau = 0}^{t-1} f_0(\x_\tau)= \frac{1}{t}\sum_{\tau = 0}^{t-1} \left(f_0(\x_\tau) + {\ip{\lam_\tau,\f(\lam_\tau)}}\right) - \frac{1}{t}\sum_{\tau = 0}^{t-1} {\ip{\lam_\tau,\f(\lam_\tau)}}\\
						& = \frac{1}{t}\sum_{\tau = 0}^{t-1} g(\lam_\tau) - \frac{1}{t}\sum_{\tau = 0}^{t-1} {\ip{\lam_\tau,\f(\lam_\tau)}} \geq g(\ls) - \frac{1}{t}\sum_{\tau = 0}^{t-1} {\ip{\lam_\tau,\f(\lam_\tau)}}. \label{fzbnd} 
						\end{align}
						\normalsize
						Next, the second term in \eqref{fzbnd} can be bounded as 
						\small
						\begin{align}
						\hspace{0mm}\norm{\lam_{t+1}}^2 &= \norm{\pll{\lamt - \mu \ft(\lamt)}}^2 \leq \norm{\lamt - \mu \ft(\lamt)}^2 \label{non_exp}\\
						& \leq \norm{\lamt}^2 \!\!\!- 2\mu {\ip{\lamt,\f(\lamt)}} \!\!+\!\! \mu^2\norm{\ft(\lamt)}^2 -2\mu{\ip{\ft(\lamt)-\f(\lamt),\lamt}}\\
						\Rightarrow 2\mu{\ip{\lamt,\f(\lamt)}} &\leq \norm{\lam_t}^2 - \norm{\lam_{t+1}}^2 + \mu^2G^2 -  2\mu{\ip{\ft(\lamt)-\f(\lamt),\lamt}}\label{zero_lam} 
						\end{align}
						\normalsize
						where, \eqref{non_exp} follows from the non-expansiveness property of the projection operator and from the fact that $\boldsymbol{0}\in\Lambda$ $\left(\norm{\lam_{t+1}}=\norm{\lam_{t+1}-\boldsymbol{0}}\right)$, and \eqref{zero_lam} follows from (\textbf{A3}). Taking sum over $\tau = 0, 1, \ldots, t-1$ and dividing by $2\epsilon t$ yields
						\small
						\begin{align}
						{\frac{1}{t} \sum_{\tau = 0}^{t-1}} {\ip{\lam_\tau,\f(\lam_\tau)}} &\leq \frac{\norm{\lam_{0}}^2}{2\mu t}-\frac{\norm{\lam_{t+1}}^2}{2\mu t}  + \frac{\mu G^2}{2} -  {\frac{1}{t}\sum_{\tau = 0}^{t-1}} {\ip{\f_\tau(\lam_\tau)-\f(\lam_\tau),\lamt}} \leq \frac{\norm{\lam_0}}{2n} + \frac{\epsilon G^2}{2} + C'_t(n,\epsilon) \label{pbnd}
						\end{align}
						\normalsize
						where, 
						%\begin{align} \label{ctne2}
					$	C'_t(n,\epsilon) := \abs{\frac{1}{ t}\sum_{\tau = 0}^{t-1} {\ip{\f_\tau(\lam_\tau)-\fa(\lam_\tau),\lam_\tau}}}.$
					%	\end{align}
						The bound in \eqref{t2} follows by plugging back \eqref{pbnd} into \eqref{fzbnd}. The analysis for $C'_t(n,\epsilon)$ is much the same as in the proof of Theorem \ref{conrate}. The only difference for this case is that the iterate bound becomes $\norm{\lam_t} \leq \sqrt{K}\lambda_{\max}$ from \eqref{dualupp}. Consequently, after rearranging various terms in $C'_t(n,\epsilon)$ and using the triangle inequality in \eqref{slbnd}, \eqref{slbnd2}, and \eqref{sl2}, all occurrences of $\Lambda_{\max}$ get replaced with $\sqrt{K}\lambda_{\max}$. Since this is equivalent to redefining the constant $\Lambda_{\max}$ appropriately, the required rate results continue to hold. 
					\end{IEEEproof}\vspace{-4mm}
							\section{Application to D2D Communications}\label{d2d}
							This section details some implementation aspects of the SDSD algorithm in the context of D2D communication problem considered in this paper under slow and fast fading scenarios. The Assumptions (\textbf{A1})-(\textbf{A4}) are also verified for the problems at hand so as to ensure that Theorems 1 and 2 hold. 
							
							Before proceeding, the {SDSD} algorithm for the general form of the D2D problem \eqref{d2d0} is detailed. Specifically, the Lagrangian is given by 
							\begin{align}\label{lagrangian}
							\hspace{-5mm}L(r,\{p^i\}_{{i\in\mathcal{M}}},\lambda)=U(r)-&\Ex{\sum\limits_{{i\in\mathcal{M}_t}}c_t^i p^i_t}+\lambda\Ex{\sum\limits_{{i\in\mathcal{M}_t}}R_i(p^i_t,\gamma^i_t)-r}
							\end{align}\normalsize
							which yields the following stochastic algorithm. \colb{Since the Lagrangian is separable in $r$ and $p^i_t$}, starting with arbitrary $\lambda_0$, the primal iterates at time slot $t$ become:
							\begin{align}
							{r}_t(\lambda_t) & \in  \argmax_{r_{\min}\leq r\leq r_{\max}}  U(r)-\lambda_t r\label{primal_update_gen}\\
							\{p_t^i(\lambda_t)\}_{{i\in\mathcal{M}_t}} &\in \argmax_{\{\mathring{p}^i\}_{{i\in\mathcal{M}_t}} \in {\Pi_t}} \sum\limits_{{i\in\mathcal{M}_t}}\left[\lambda_t R_i(\mathring{p}^i,\gamma^i_t)-c_t^i\mathring{p}^i\right]. \label{power_alloc}
							\end{align}
							At the end of each time slot, the dual variable is updated as
							\begin{align}\label{dual_update_gen}
							\lambda_{t+1}=\pl{\lambda_{t}-\epsilon \left[\sum_{{i\in\mathcal{M}_t}} R_{i}(p^{i}_t(\lambda_t),\gamma^{i}_t)-{r}_t(\lambda_t)\right]}
							\end{align}
							Recall that the set of functions $\mathcal{P}$ is such that only one user, denoted by 
							\vspace{-5mm}
							\colb{\begin{align}\label{user_selection}
							i_t:=\arg\max_{i\in\mathcal{M}_t} p^i_t(\lambda_t),
							\end{align}} is allocated non-zero power at time slot $t$. Therefore, the dual variable is updated as
			\vspace{-6mm}
							\begin{align}\label{dual_update_gen2}
							\hspace{-0.2cm}\lambda_{t+1}&=\pl{\lambda_{t}-\epsilon \left[R_{i_t}(p^{i_t}_t(\lambda_t),\gamma^{i_t}_t)-{r}_t(\lambda_t)\right]}
							\end{align}
							The full algorithm is summarized in Algorithm \ref{algo1}. 
							\begin{algorithm} \label{Algo1}
								\caption{Edge Caching via D2D Communications}\label{algo1}
								\begin{algorithmic}[1]
									\STATE {\bf Initialize} $\lambda_0$ \\
									\textbf{Repeat for $t ={ 0, 1}, 2, \ldots, $}
									\STATE \hspace{0.5cm} Collect $c^i_t$ and $\gamma^i_t$ from each active UE ${i\in\mathcal{M}_t}$
									\STATE \hspace{0.5cm} \textbf{Primal update} Determine the ``winning'' cache $i_t$ and the allocated power \colb{from \eqref{power_alloc}} 
									\STATE \hspace{0.5cm} Download at rate $R_{i_t}(p^{i_t}_t(\lambda_t),\gamma^{i_t}_t)$ from user $i_t$
									\STATE \hspace{0.5cm} \textbf{Dual update} Update $\lambda_{t+1}$ using \eqref{dual_update_gen2}
								\end{algorithmic}
							\end{algorithm} 							
							The rate analysis developed in Sec. \ref{conv} applies to the present problem under the following assumptions.
							\begin{enumerate}
								\item[\textbf{B1}. ] \emph{Continuous random variables: } The random variables $\bg_t = ({\mathcal{M}_t}, \{c^i_t\}_{{i\in\mathcal{M}_t}}, \{\gamma^i_t\}_{{i\in\mathcal{M}_t}})$ are i.i.d., have continuous cdfs, and finite supports, i.e., ${\mathcal{M}_t} \subset \mathcal{M}$, $\gamma_t^i \in [\gamma_{\min}, \gamma_{\max}]$, and $c^i_t \in [c_{\min} , c_{\max}]$ for each ${i\in\mathcal{M}_t}$.
								\item[\textbf{B2}. ] \emph{Power constraints:} The set $\mathcal{P} := \{\p:\Rn^{3M}\rightarrow \Rn^M \mid \p_{\bg} \in \Pi_{\bg}\}$, where for any $\bg_t$, we have that $\Pi_{\bg_t} := {\Pi_t}=\{\p_{t} \in \Rn^{M} \mid p_t^j = 0 ~ {j \notin {\mathcal{M}_t}}, \norm{\p_{t}}_0 = 1, p^{i_t}_tc_t^{i_t} \in [C_{\min}, C_{\max}] \}$,
								where $i_t := \arg\max_{{i\in\mathcal{M}_t}} p^i_{t}$ and $p^i_t = [\p_t]_i$. 
								\item[\textbf{B3}. ] \emph{High SNR: } It is assumed that $\gamma^i_t \gg 1$ for all ${i\in\mathcal{M}_t}$. 
							\end{enumerate}
							
							The finite support of the random quantities is again motivated from practical considerations. The set $\mathcal{P}$ also includes limits on the maximum affordable cost $C_{\max}$ and the minimum operational cost or minimum allowable transaction amount $C_{\min}$. A maximum power constraint of the form $p^i_t \leq P_{\max}$ may also be included within $\mathcal{P}$. However, for the present application, it is assumed that the caches are not energy constrained, so that $P_{\max} \gg C_{\max}/c^i_t$ for all ${i\in\mathcal{M}_t}$. In other words, the user's cost constraint is much more stringent than the cache's energy constraint. Finally, the high SNR assumption is justified if there are always enough mobile caches available at all slots. In a typical setting, the {MoUE} may ``see'' hundreds of advertisements from potential mobile cache servers, but may choose to consider only tens of users with which control messages may be exchanged easily. Next, the discussion for slow and fast fading cases will be carried out.							
							\subsection{Slow Fading}\label{slow}
							Recall that under slow fading, since power allocation occurs every coherence interval, \colb{we have for high SNR (cf. (\textbf{B3})),} that $R_i(p_t^i(\lambda_t),\gamma_t^i) :\approx W\colb{\log_2} (p_t^i(\lambda_t)\gamma_t^i/\alpha)$. The primal iterate in \eqref{power_alloc} can be found in two steps. First the optimum transmit power for all potential users is determined, i.e., for each ${i\in\mathcal{M}_t}$, \vspace{-5mm}
							\begin{align}
							\hat{p}_t^i(\lambda_t) &=\argmax_{\mathring{p}^i} ~ \lambda_t  W\colb{\log_2} (\mathring{p}^i\gamma_t^i/\alpha)-c_t^i\mathring{p}^i \\
							& \hspace{1cm} \text{s. t. } ~~~ C_{\min}\leq c_t^i\mathring{p}^i\leq C_{\max} \nonumber\\
							&= \left[\frac{W\lambda_t}{c_t^i}\right]_{C_{\min}/c_t^i}^{C_{\max}/c_t^i} \label{primal_slow2}
							\end{align}
							The winning user is the one that maximizes the objective function, i.e., 
							\begin{align}
							\hspace{-0.2cm}	i_t &=\argmax_{{i\in\mathcal{M}_t}}  [\lambda_t R_i(\hat{p}^i_t,\gamma^i_t)-c_t^i\hat{p}^i_t]= \argmax_{{i\in\mathcal{M}_t}} \frac{\gamma^i_t}{c^i_t} \label{itlam1}
							\end{align}
							where the expression in \eqref{itlam1} derived in Appendix \ref{itlam}. An implication of \eqref{itlam1} is that the random variable $i_t$ is i.i.d. Finally, it holds that $p^j_t = \hat{p}^{i_t}_t(\lambda_t)$ for $j = i_t$ and zero otherwise. Similarly,  $r_t$ is calculated as
							%\begin{align}\label{primal_slow}
						$	{r}_t(\lambda_t) =\left[\frac{1}{\lambda_t}\right]_{r_{\min}}^{r_{\max}}$
						%	\end{align}
							resulting in the dual update
							\begin{align}\label{dual_slow}
							\lambda_{t+1}=\pl{\lambda_{t}-\epsilon \left[W\colb{\log_2}(p^{i_t}_t(\lambda_t)\gamma^{i_t}_t/\alpha)-{r}_t(\lambda_t)\right]}.
							\end{align}
							
							An additional assumption regarding the parameter values is made in the slow fading case:
							
							\noindent \textbf{B4}. ~ \emph{Strict feasibility: } The problem parameters satisfy $r_{\min} < \Ex{\max_i \colb{\log_2}(C_{\max}\gamma_t^i/c^i_t)}$. 
							
							The strict feasibility condition is required for ensuring the existence of a Slater point. Since it holds that ${\gamma_t^i/c^i_t \geq {\gamma_{\min}/c_{\max}}}$, it is possible to satisfy (\textbf{B4}) by keeping $r_{\min}$ sufficiently small and/or if $\gamma_{\min}$ is sufficiently large. 
							
							Having stated the algorithm and all required assumptions, the following Lemma summarizes the main result of this subsection. 							
							\begin{cor} \label{lemslow}
								Under (\textbf{B1})-(\textbf{B4}), the iterates obtained from \eqref{primal_slow2}-\eqref{dual_slow} adhere to the rate bounds stated in Theorems 1 and 2.
							\end{cor}
							\begin{IEEEproof}
								For the results in Theorem 1 and 2 to apply, it suffices to verify that assumptions (\textbf{A1})-(\textbf{A4}) are satisfied under the slow fading case. The random variable $\bg_t$ has a non-atomic pdf since the channel gains $\gamma^i_t$ have a continuous cdf, thus confirming (\textbf{A1}); see also \cite{rao1972remark}. The Slater's condition is met by choosing $\tilde{r} = r_{\min}$ and $\tilde{p}^{i_t}_t = C_{\max}/c_t^{i_t}$ where $i_t$ is given in \eqref{itlam1} and zero for all $ j \neq i_t$. For such a choice, it holds from (\textbf{B4}) that $\tilde{r} < \Ex{\colb{\log_2}(\tilde{p}^{i_t}_t\gamma^{i_t}_t)}$, \colb{which is the} required condition for strict feasibility. For a given $\lambda$, the subgradient function is given by 
								\begin{align}
								f_t(\lambda) &= W\colb{\log_2}({p}^{i_t}_t(\lambda)\gamma^{i_t}_t/\alpha)-r_t(\lambda)
								\end{align}
								where $i_t$ and $p^{i_t}_t$ are evaluated as in \eqref{itlam1} and \eqref{primal_slow2}. A 
								bound on the subgradient (cf. (\textbf{A3})) may therefore be found as
$\abs{f_t(\lambda)}  \leq W\colb{\log_2}\left(\frac{C_{\max}\gamma_{\max}}{\alpha c_{\min}}\right) + r_{\max} =: G$. Next, {in order to verify (\textbf{A4})}, the expression for the stochastic subgradient error $e_t(\lambda):=f_t(\lambda)-\Ex{f_t(\lambda)}$ is first derived. Recalling that $i_t = \argmax_i \gamma^i_t/c^i_t$, consider the following three cases, 
								\begin{enumerate}
									\item When $\lambda < C_{\min}/W$, it holds that $p^{i_t}_t = C_{\min}/c^{i_t}_t$, implying that
									\small
									\begin{align}
									e_t(\lambda) &= W \colb{\log_2}\left(\frac{C_{\min}\gamma^{i_t}_t}{\alpha c^{i_t}_t}\right) - \left[\frac{1}{\lambda}\right]_{{r_{\min}}}^{r_{\max}} - \Ex{W \colb{\log_2}\left(\frac{C_{\min}\gamma^{i_t}_t}{\alpha c^{i_t}_t}\right) - \left[\frac{1}{\lambda}\right]_{{r_{\min}}}^{r_{\max}}} \\
									&= W \colb{\log_2}\left(\frac{\gamma^{i_t}_t}{c^{i_t}_t}\right) - \Ex{W \colb{\log_2}\left(\frac{\gamma^{i_t}_t}{c^{i_t}_t}\right)}.
									\end{align} \normalsize
									where the expectations are with respect to $\bg_t$. 
									\item When $C_{\min}/W \leq \lambda \leq C_{\max}/W$, it holds that $p^{i_t}_t = W\lambda\gamma^{i_t}_t/c^{i_t}_t$, implying that
									\small
									\begin{align}
									e_t(\lambda) &= W \colb{\log_2}\left(\frac{W\lambda\gamma^{i_t}_t}{c^{i_t}_t}\right) - \left[\frac{1}{\lambda}\right]_{{r_{\min}}}^{r_{\max}}- \Ex{W \colb{\log_2}\left(\frac{W\lambda\gamma^{i_t}_t}{c^{i_t}_t}\right) - \left[\frac{1}{\lambda}\right]_{{r_{\min}}}^{r_{\max}}} \\
									&= W \colb{\log_2}\left(\frac{\gamma^{i_t}_t}{c^{i_t}_t}\right) - \Ex{W \colb{\log_2}\left(\frac{\gamma^{i_t}_t}{c^{i_t}_t}\right)}.
									\end{align} 
									\normalsize
									\item Similarly, when $\lambda > C_{\max}/W$, it holds that $p^{i_t}_t = C_{\max}/c^{i_t}_t$, implying that
									\begin{align}\label{final_diff}
									e_t(\lambda) &= W \colb{\log_2}\left(\frac{\gamma^{i_t}_t}{c^{i_t}_t}\right) - \Ex{W \colb{\log_2}\left(\frac{\gamma^{i_t}_t}{c^{i_t}_t}\right)}.
									\end{align} 
								\end{enumerate}
								Therefore, the subgradient error is a zero-mean random variable that does not depend on $\lambda$, and is therefore trivially continuously differentiable in $\lambda$. 
							\end{IEEEproof}
							\vspace{-5mm}
							\subsection{Fast Fading}\label{fast}
							In the more realistic fast fading case, the power allocation and downloads occur over several coherence intervals.  Under the high SNR assumption, the rate becomes , $R_i(p^i_t,\gamma^i_t) \approx W\colb{\log_2}(p_t^i\gamma_t^i/\alpha)+W\psi_i$ where $\psi_i = \Eh{\colb{\log_2}(h_i)}$ for a given user $i$ \cite{tse}. As in the slow fading case, the primal iterates are again found in two steps. First, the power allocation for a potential user $i$ is found,
							\begin{align}
							\hat{p}_t^i(\lambda_t) &= \left[\frac{W\lambda_t}{c_t^i}\right]_{C_{\min}/c_t^i}^{C_{\max}/c_t^i}. \label{primal_fast2}
							\end{align}
							It is shown in Appendix \ref{itlam} that the winning user for the fast fading case can be written as
							\begin{align}
							i_t &= \argmax_{i\in {\mathcal{M}_t}} \colb{\log_2}\left(\frac{\gamma^i_t}{c^i_t}\right) + \psi_i. \label{itlam2}
							\end{align}
							Finally, since $r_t(\lambda_t) = {\max\{\min\{1/\lambda_t, r_{\max}\},r_{\min}\}}$ as before, the dual update is given by
							\begin{align}\label{dual_fast}
							\lambda_{t+1}\!\!=\!\pl{\!\lambda_{t}\!-\!\epsilon \left[\!W\!\colb{\log_2}(p^{i_t}_t(\!\lambda_t\!)\gamma^{i_t}_t\!/\alpha)\! +\! W\psi_{i_t}\!\!-{r}_t(\!\lambda_t\!)\!\right]\!}.
							\end{align}
							In order to apply the rate bounds in Theorems 1 and 2, we again assume (\textbf{B1})-(\textbf{B3}), and make the following assumption analogous to (\textbf{B4}). 
							
							\noindent \textbf{B5}. ~ \emph{Strict feasibility: } The problem parameters satisfy $r_{\min} < \Ex{\max_i \left\{\colb{\log_2}(C_{\max}\gamma_t^i/c^i_t) + \psi_i\right\}}$. 
							
							As in the slow fading case, (\textbf{B5}) allows us to obtain a Slater point, as required by (\textbf{A2}). The following Lemma summarizes the result for the fast fading case. 
							
							\begin{cor} \label{lemfast}
								Under (\textbf{B1})-(\textbf{B3}) and (\textbf{B5}), the iterates obtained from \eqref{primal_fast2}-\eqref{dual_fast} adhere to the rate bounds stated in Theorems 1 and 2. 
							\end{cor}
							\begin{IEEEproof}
								The 
								As in Lemma \ref{lemslow}, it suffices to verify assumptions (\textbf{A1})-(\textbf{A4}). The random variable $\bg_t$ has a non-atomic pdf as remarked earlier. Similarly, it can be verified that the Slater point is given by $\tilde{r} = r_{\min}$ and $\tilde{p}^{i_t}_t = C_{\max}/c_t^{i_t}$ where $i_t$ is as given in \eqref{itlam2}, and zero for all $ j \neq i_t$. The subgradient bound required in (\textbf{A3}) now becomes, $\abs{f_t(\lambda)} \leq W\colb{\log_2}\left(\frac{C_{\max}\gamma_{\max}}{\alpha c_{\min}}\right) + W\psi_{\max} + r_{\max} =: G$ where $\psi_{\max} := \max_{i} \psi_i$. Finally, in order to verify (\textbf{A4}), we proceed as in the proof of Lemma \ref{lemslow} and derive an expression for the subgradient error $e_t(\lambda):=f_t(\lambda)-\Ex{f_t(\lambda)}$. Since expression for the allocated power is the same for the two cases, it can be seen that for the fast fading case as well $e_t(\lambda) = W \colb{\log_2}\left(\frac{\gamma^{i_t}_t}{c^{i_t}_t}\right) - \Ex{W \colb{\log_2}\left(\frac{\gamma^{i_t}_t}{c^{i_t}_t}\right)}$ where $i_t$ is found as in \eqref{itlam2}. Since $e_t(\lambda)$ does not depend on $\lambda$, (\textbf{A4}) also holds trivially in the fast fading scenario.  
							\end{IEEEproof}
								\begin{figure*}
									\setcounter{subfigure}{0}
									\begin{subfigure}{0.5\columnwidth}
										\includegraphics[width=\linewidth, height = 0.5\linewidth]
										{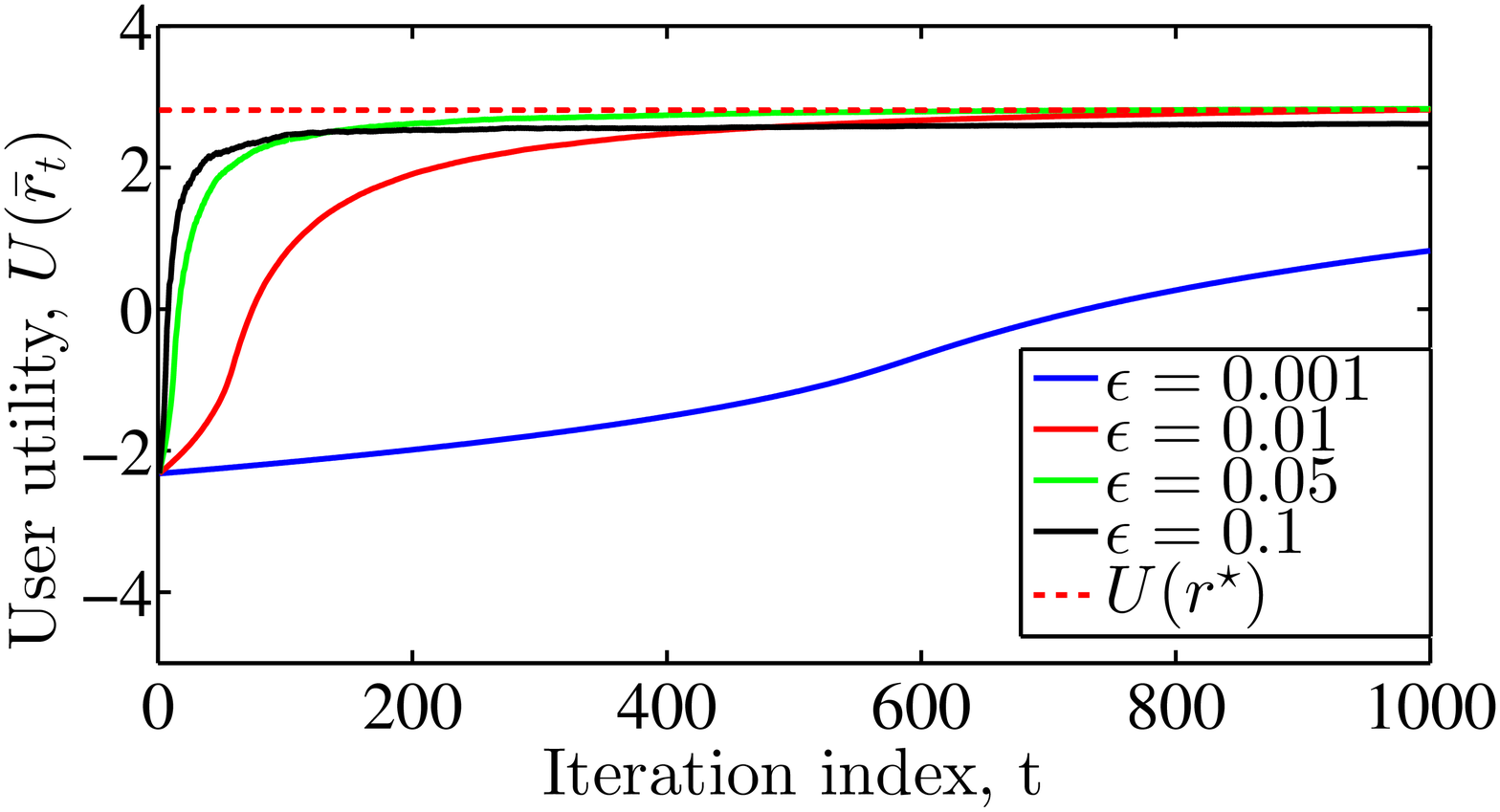}
										\caption{The user utility}
										\label{fig1_top}
									\end{subfigure}
									\begin{subfigure}{0.5\columnwidth}
										\includegraphics[width=\linewidth, height = 0.5\linewidth]
										{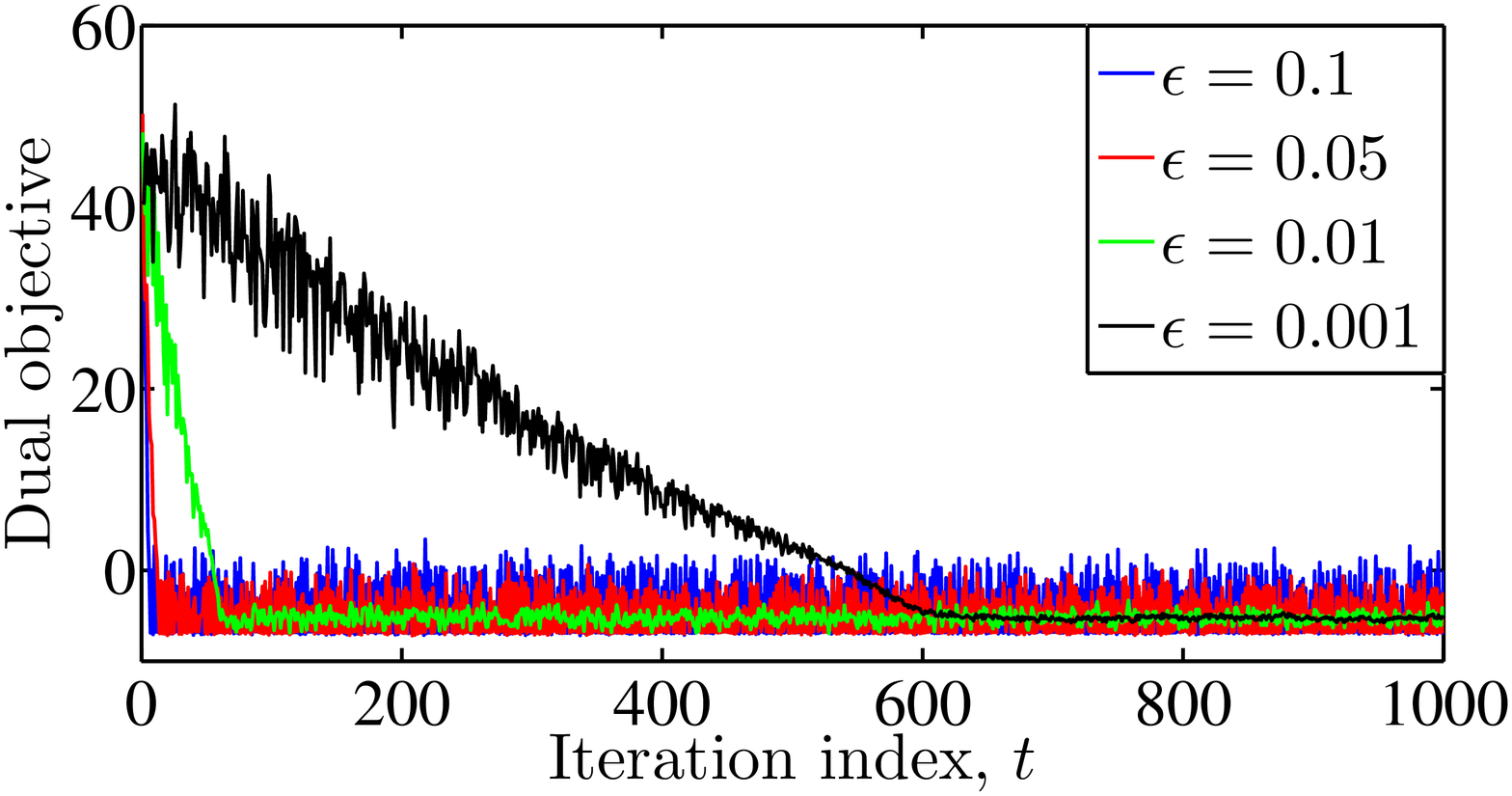}
										\caption{The dual objective}
										\label{fig1_bottom}
									\end{subfigure}\vspace{-4mm}
										\caption{Performance behavior} \vspace{-12mm}
								\end{figure*}
								\vspace{-5mm}
							\section{Numerical Tests}\label{num}
							This section describes the numerical tests on the D2D example discussed in Sec. \ref{d2d}. The convergence rate of the SDSD algorithm is studied for the fast fading scenario depicted in Fig. 2. {For the simulations, we consider $M\!=\!25$ operational UEs.} At each time slot, the MoUE receives advertisements from a random subset $\mathcal{M}_t$ of 5 to 25 UEs. Without loss of generality, downloading from the $i$-th UE incurs a cost of {$c^i_t\! =\! i$} per unit of transmit power. The lower and upper limits for each transaction are set as {$C_{\min}= 1$ and $C_{\max}= 25$}, respectively. The average channel gains $\gamma^i_t$ are assumed to be Rayleigh distributed {with {$\gamma_{\min}=0.1$} and $\gamma_{\max}=65$}, and for simplicity, the parameters $\alpha$, and $\psi_i$ are all set to unity. In order to keep the numbers realistic, we set {$W=1$ MHz}. Finally, in order to ensure Slater's condition, we set {$r_{\min} = 0.2$ and $r_{\max} = 10$}. {In realistic scenarios, since the optimal rate is expected to be greater than $r_{\min}$, it follows from the definition of $r_t(\lambda_t)$ in Sec.V-A  that $\lambda^\star > 1/r_{\min}$. Therefore it is safe to take  $\lambda_{\max} \gg 1/r_{\min}$.}

%							\begin{figure*}
%								\centering
%								\begin{minipage}{.5\textwidth}
%									%\centering
%									\includegraphics[width = 1\columnwidth,height=0.25\textheight]{result1_combined_2}
%									\captionsetup{font=scriptsize}
%									\caption{The user utility $U(\bar{r}_t)$ (top), Dual objective (bottom).}
%									\label{fig1}
%								%	\vspace{-5mm}
%								\end{minipage}%
%%								\hspace{2mm}\begin{minipage}{.45\textwidth}
%%								%	\centering
%%									\includegraphics[width = \columnwidth,height=0.2\textheight]{result4_3_1}
%%																	\captionsetup{font=scriptsize}
%%																		\caption{Behavior of $C_t(n,\epsilon)$ for fixed $n$ (left) and fixed $\epsilon$ (right) }
%%									\label{fig4_1}
%%								\end{minipage}
%																\vspace{-0mm}	
%							\end{figure*}

%							\begin{figure}
%								\centering
%								\includegraphics[width = 0.5\columnwidth,height=0.25\textheight]{result1_combined_2}
%								\caption{The user utility $U(\bar{r}_t)$ (top), Dual objective (bottom).}
%								\label{fig1}
%								\vspace{-5mm}
%							\end{figure}							
							Fig. \ref{fig1_top} shows the evolution of the utility function calculated using running averages  $\bar{r}_t := 1/t\sum_{\tau = 0}^{t-1} r_\tau$ of the allocated rate with iterations. As expected from Theorem \ref{conv_primal}, the utility function converges to a value that is closer to the optimal when $\epsilon$ is small. Similarly, Fig. \ref{fig1_bottom} shows the evolution of the dual objective function, which again converges to a point closer to the optimal when $\epsilon$ is small. Observe from the results that for $\epsilon = 0.1$, the oscillations continue even as number of iterations go to infinity, as implied by Theorem \ref{conrate}. {These oscillations are allowed due to the presence of an $\O(\eps)$ term on the right-hand side of (17), and are well-documented for the constant step size stochastic subgradient type algorithms \cite{bottou1998online,nedic2001convergence,ale10}.}

\begin{figure}
									\setcounter{subfigure}{0}
									\begin{subfigure}{0.5\columnwidth}
										\includegraphics[width=\linewidth, height = 0.5\linewidth]
										{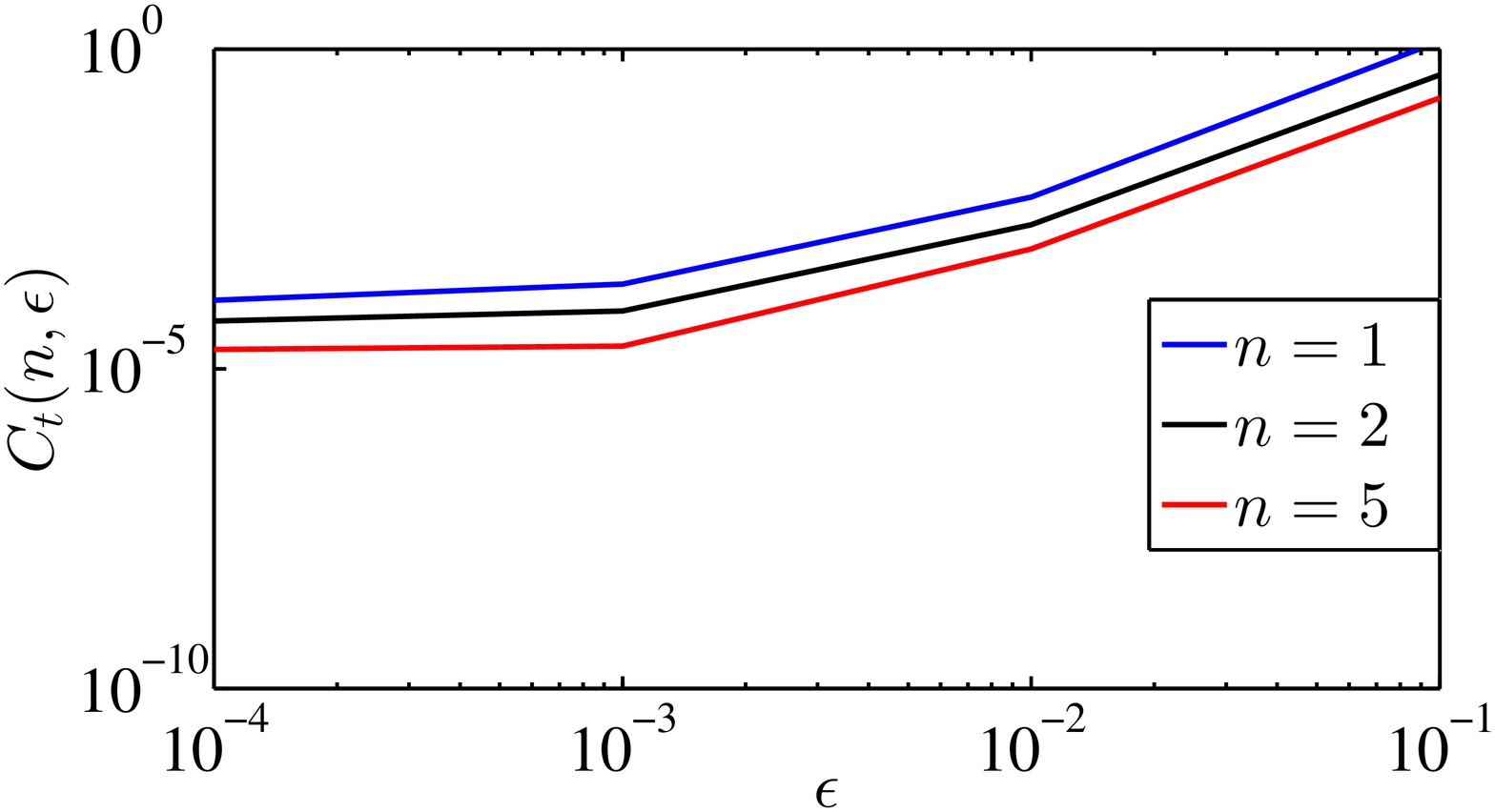}
										%\caption{Fixed $n$ }
										%\label{fig00}
									\end{subfigure}
									\begin{subfigure}{0.5\columnwidth}
										\includegraphics[width=\linewidth, height = 0.5\linewidth]
										{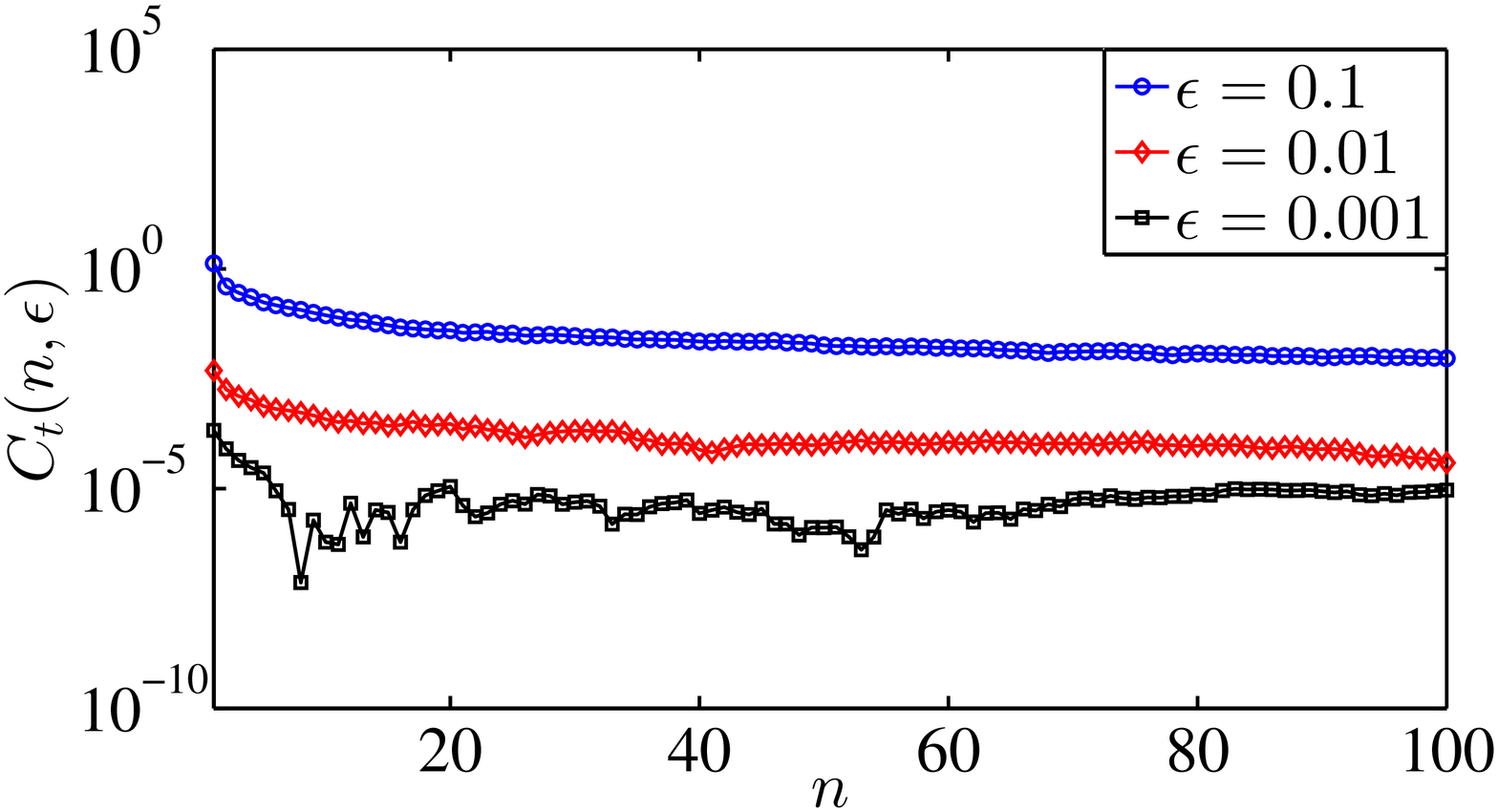}
										%\caption{Fixed  $\epsilon$}
										%\label{fig11}
									\end{subfigure}
									\vspace{-5mm}
									\caption{Behavior of $C_t(n,\epsilon)$ for (a) fixed $n$ and (b) fixed $\epsilon$}\vspace{-8mm}\label{fignep}
								\end{figure}

							The convergence rate result of Theorem \ref{conrate} is further illustrated in \colb{Fig. \ref{fignep}(a) and Fig. \ref{fignep}(b)}, where the deterministic terms in \eqref{conrateeq} are not included. The stochastic term $C_t(n,\epsilon)$ is calculated from {\eqref{ctne}} and plotted against both $\epsilon$ and $n$. It can be seen from both plots that $C_t(n,\epsilon) \rightarrow 0$ as either $n \rightarrow \infty$ or $\epsilon \rightarrow 0$, as claimed in Theorem \ref{conrate}. 
							\begin{table}[h]
								\centering
								\vspace{-4mm}
								\begin{tabular}{ |c|c|c|c| } 
									\hline
									UE selection & Downloaded data (Mb) & Cost incurred &  Avg. Utility minus penalty  \\
									\hline
									Proposed &    16316 & 4501 &   0.79\\ 
									\hline
									Random  &  12465 & 19951&    -3.20\\
									\hline
									Opportunistic &   3987 & 177 & 0.67\\
									\hline
								\end{tabular}
								\vspace{-1mm}
								\caption{\colb{Performance comparison ($1000$ times slots)}} \label{table}
							\end{table}
							\vspace{-8mm}
							
\colb{Having studied the convergence properties of the SDSD algorithm, we now focus on some of the nuances of the edge-caching formulation in \eqref{d2d0}. To begin with, the performance of the proposed scheme is compared against that obtained from two naive algorithms: random and opportunistic. The maximum transmit power for the three cases is scaled so as to ensure equal aggregate power consumption. As the name suggests, an MoUE following the \emph{random} scheme selects an available UE randomly and without paying any attention to the channel or the cost of the UE. The data is transmitted at the maximum power so as to ensure the maximum rate. As evident from Table \ref{table}, such a scheme is able to obtain a higher download rate but also at a significantly higher cost. In contrast, the \emph{opportunistic} scheme advocates a parsimonious approach wherein the MoUE always selects an available UE with the lowest cost. Subsequently, the UE transmits with the maximum power but ultimately achieves a lower aggregate download rate, due to suboptimal channel conditions.}

\begin{figure}
										\setcounter{subfigure}{0}
										\begin{subfigure}{0.5\columnwidth}
												\includegraphics[width = \columnwidth,height=0.2\textheight]{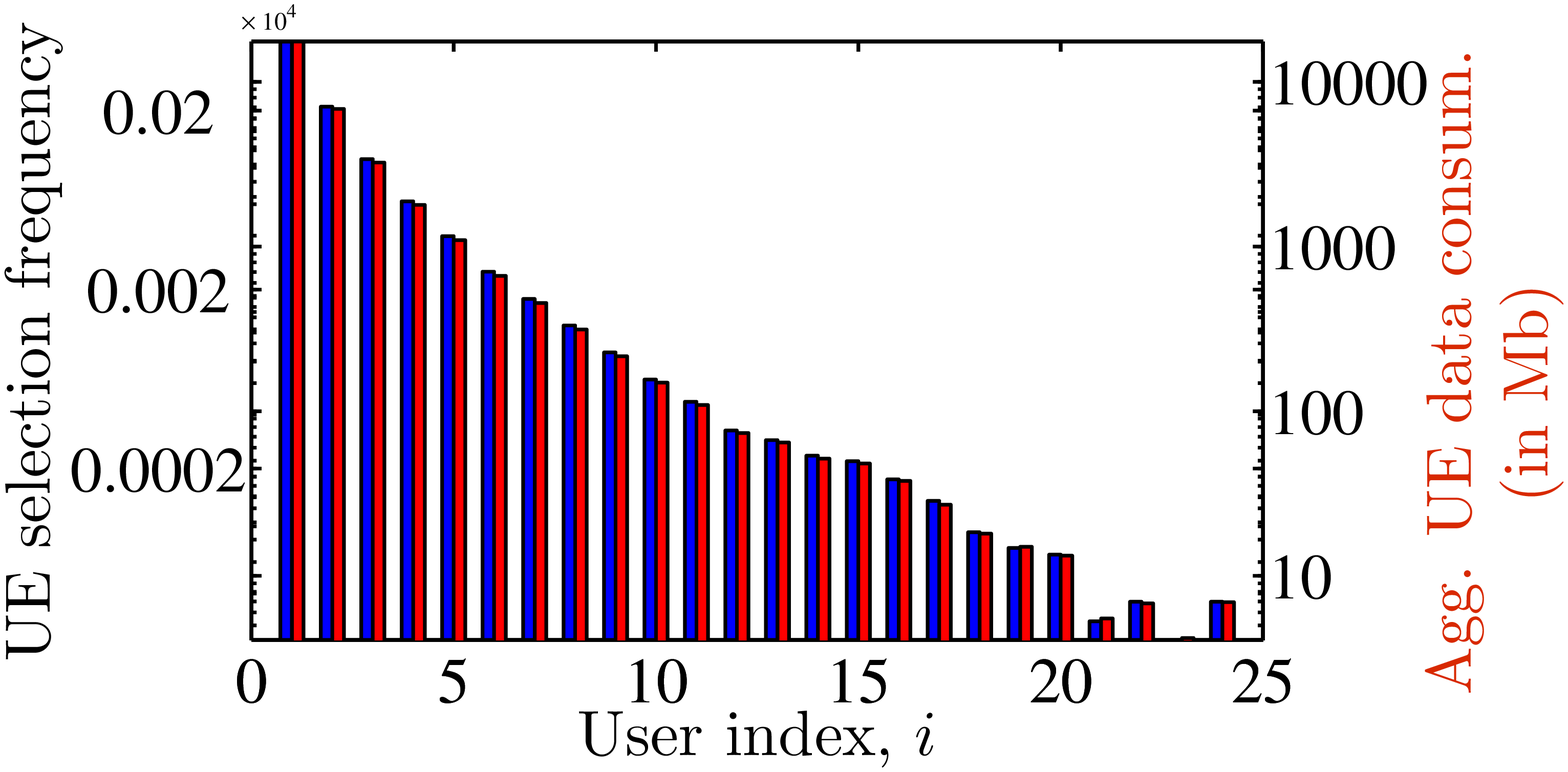}
												\captionsetup{font=scriptsize}
											%	\caption{\colb{User selection behavior and aggregate data consumption}}
												%\label{histogram}
										\end{subfigure}
										\begin{subfigure}{0.5\columnwidth}
										\includegraphics[width = \columnwidth,height=0.18\textheight]{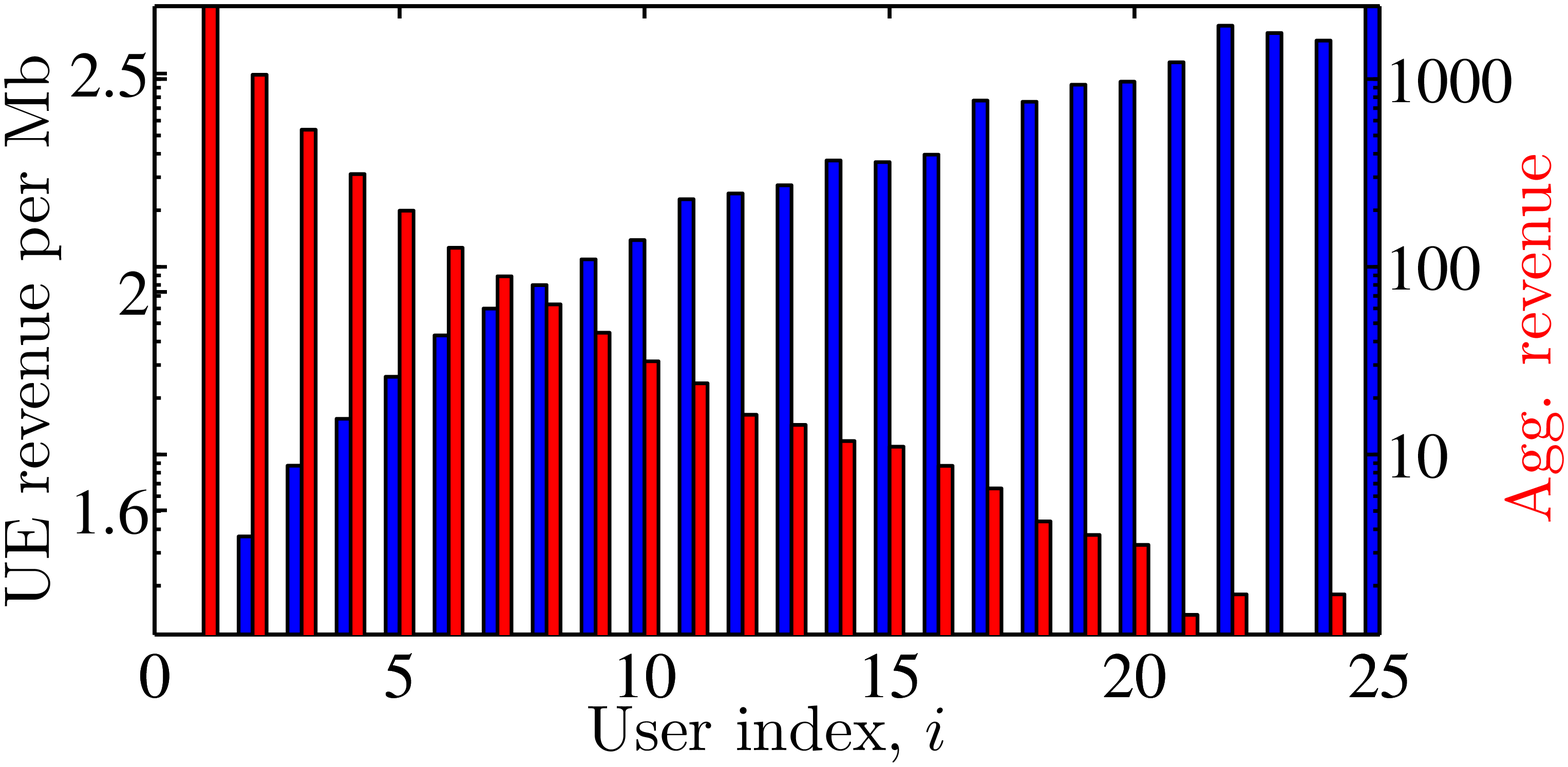}
										\captionsetup{font=scriptsize}
										%\caption{\colb{Generated revenue of each UE}}
										%\label{cost}
										\end{subfigure}\vspace{-5mm}		
%											\begin{subfigure}{0.5\columnwidth}
%												\includegraphics[width = \columnwidth,height=0.2\textheight]{result7}
%												\captionsetup{font=scriptsize}
%												\caption{}
%												\label{comparison_rate}
%											\end{subfigure}
%												\begin{subfigure}{0.5\columnwidth}
%													\includegraphics[width = \columnwidth,height=0.2\textheight]{result8}
%													\captionsetup{font=scriptsize}
%													\caption{}
%													\label{comparison_cost}
%												\end{subfigure}
										\caption{\colb{(a) User selection behavior and aggregate data consumption (b) Generated revenue of each UE.}}\vspace{-8mm}\label{ue}
									\end{figure}

\colb{Fig. \ref{ue} provides results from the perspective of the UEs and is generated by running the same algorithm for 1000 independent identically distributed MoUEs. In particular, if an MoUE follows the optimal policy determined by \eqref{d2d0}, the UEs may be interested in knowing a reasonable price value. As expected, it is clear from Fig. \ref{ue}, that the UEs that charge more are selected less often have lower data usage. Consequently, the aggregate revenue of the UEs with the lowest charges is also the highest. More interestingly however, such high-priced UEs have a very high revenue per Mb of data served. The intuition here is that UEs with high costs are only selected when their channel gains are proportionally higher than the others. Therefore, all transmissions to such UEs occur at higher rates and correspondingly lower power. In summary, by operating only under favorable channel conditions, the high-priced UEs extract a greater revenue for every bit that they serve. Note however that the revenue appears to saturate, and increases very slowly for very high prices. }						

%							\begin{figure}
%								\centering
%								\includegraphics[width = 0.5\columnwidth,height=0.15\textheight]{result4_3_1}
%								\caption{Behavior of $C_t(n,\epsilon)$ for fixed $n$ (left) and fixed $\epsilon$(right) }
%								\label{fig4_1}
%								\vspace{-5mm}
%							\end{figure}
%							
							\vspace{-8mm}

								\section{Conclusions}\label{conc}
							This paper considers a general stochastic resource allocation problem and solved using constant step-size stochastic subgradient descent algorithm in an online manner. A stochastic bound on the gap between the objective function and the optimum is developed and analyzed in an almost sure sense,
							generalizing the existing results. The bounds characterize the precise manner in which the optimality gap behaves for fixed and arbitrarily small step-sizes. The convergence rate analysis is also extended to a class of stochastic resource allocation problems that utilize stochastic dual subgradient descent (SDSD) iterations. Existing results on near-optimality of the primal average objective function are again generalized for convergence rate analysis. As an example, a resource allocation problem is formulated in the context of mobile caching in device-to-device communications, and solved via {SDSD}. The regularity conditions required for the rate analysis are verified, and numerical tests are provided, further substantiating the convergence rate results. 
							\vspace{-12mm}
							\appendices
							\section{A bound on $\norm{\ls}_{\infty}$} \label{lambda}
						From (\textbf{A2}), there exists some $\tilde{\x} \in \mathcal{X}$ and $\tilde{\p} \in \Pc$, such that $\Ex{\st(\tilde{\x}, {\tilde{\p}_t})} > 0$, {where recall that $\tilde{\p}_t:=\tilde{\p}_{\bg_t}$} for all ${t\in \Nn_0}$ and the expectation is with respect to $\bg_t$. Given $\lt \in \Rn^{K}_{+}$, define the sublevel set $\Q_{\lt}:=\{\lam \in \Rn^K_{+} \mid g(\lam)\leq g(\lt)\}$, and observe that for any $\lam \in \Q_{\lt}$, it holds that
						\begin{align}
						g(\lt) &\geq g(\lam)=\max_{\x \in \Xc, \p \in \Pc} f_0(\x) + {\ip{\lam,\Ex{\st({\x,\p_t})}}}\geq f_0(\tilde{\x}) + {\ip{\lam,\Ex{\st({\tilde{\x},\tilde{\p}_t})}}}.\label{bound2}
						\end{align}
						Rearranging the expression in \eqref{bound2}, we obtain 
						\begin{align}
						\sum\limits_{k=1}^{K}[\lam]_k\Ex{\left[\st({\tilde{\x},\tilde{\p}_t})\right]_k}&\leq g(\lt)-f_0(\tilde{\x})
							\Rightarrow 	\norm{\lam}_{\infty} \leq \sum\limits_{k=1}^{K}[\lam]_k\leq \frac{g(\lt)-f_0(\bar{\x})}{{\chi(\tilde{\x},\tilde{\p})}} \label{lambound2}
						\end{align}
						where ${\chi(\tilde{\x},\tilde{\p})}:=\min_{1\leq k\leq K}\Ex{\left[\st({\tilde{\x},\tilde{\p}_t})\right]_k}$. Observe that $\Q_{\ls}=\{\lam \in \Rn^K_{+} \mid g(\lam) \leq \mathsf{D}\}$, so that it follows from \eqref{lambound2} that\vspace{-5mm}
						\begin{align}\label{lsbound1}
						\norm{\ls}_{\infty} \leq \frac{\mathsf{D}-f_0(\bar{\x})}{{\chi(\tilde{\x},\tilde{\p})}}
						\end{align}
						Finally, since $g(\lt) \geq \mathsf{D}$ for all $\lt \in \Rn^K_{+}$, the bound in \eqref{lsbound1} can be relaxed to yield \eqref{lambound}.
							\vspace{-5mm}
\section{Asymptotic properties of $C_T$} \label{appbound}
\begin{IEEEproof}
In order to study the convergence rate of $C_T$, the time is divided into epochs of duration $1/\epsilon$, so that there exists some $n\! \geq \!0$ that satisfies $n/\epsilon \!\leq\! T \!<\! (n+1)/\epsilon$. Since $n:= \lfloor \epsilon t \rfloor$, {where $\lfloor\cdot\rfloor$ denotes the floor operation,} is an arbitrary number, such a split allows the value of $t$ to be increased by keeping either $n$ or $\epsilon$ fixed, and varying the other. It is therefore possible to separately study the effects of choosing larger $n$ or smaller $\epsilon$ values. For instance, if $\epsilon = 0.1$, the time is divided into epochs of duration 10 iterations each. Hence, the zeroth epoch consists of iterations $0 \leq t \leq 9$,  the first epoch consists of iterations $10 \leq t < 19$, and so on. With such a split, the classical asymptotic analysis for $t \rightarrow \infty$ is equivalent to fixing $\epsilon$ and letting $n \rightarrow \infty$. Additionally, the proposed split allows us to study the case when $n$ is fixed, but the algorithm is run with different values of $\epsilon$.							
							
								This proof is devoted to the analysis of  $C_t(n,\epsilon)$, and relies on rearranging the terms in \eqref{ctne} so that the results in developed in Lemma \ref{errbounds} can be applied. The proof is divided into two parts, one for each mode of convergence in \eqref{ctlim}.
								
								\emph{Fixed $n\! < \!\infty$ and $\mu \rightarrow 0$:}
								For this case, $C_t(n,\epsilon)$ is split into summands corresponding to each epoch till time $t$, that is,\vspace{-2mm}
								\begin{align}\label{cmmu} 
								C_t(n,\epsilon) = \abs{\frac{1}{\mu t}\sum_{m=0}^{n} C^m(\mu)}\ \ \ \text{where,}\ \ \ \ \ C^m(\mu) :=  \mu\sum_{\tau = \ell_m}^{u_m}{\ip{\f_\tau(\lam_\tau)-\f(\lam_\tau),\lam_\tau-\ls}}.
								\end{align}
								The limits in the summation are defined as $\ell_m := m/\mu$ and $u_m := \frac{(m+1)}{\mu}-1$ for $m < n$ while $u_n:=t-1$. Next, define for all $\lam \in\Lambda$ and $\ell_m  \leq \tau \leq u_m$ 
								\begin{align}\label{ztau}
								z_{\tau}(\lam) = \mu\sum_{\iota = \ell_m}^{\tau} {\ip{\f_\iota(\lam) - \f(\lam),\lam-\ls}}.
								\end{align}
								Substituting \eqref{ztau} in \eqref{cmmu}, we obtain 
								\begin{align}\label{}
								C^m(\mu) = z_{u_m}(\lam_{u_m+1}) - \sum_{\tau = \ell_m}^{u_m} \left(z_\tau(\lam_{\tau+1}) - z_\tau(\lam_\tau)\right).
								\end{align}
								Such a split allows us to use \eqref{a5a6} in order to bound the magnitude of each term separately. Specifically, letting $\T_m := \{\ell_m ,\ell_m+1, \ldots, u_m\}$, 
								\begin{align}
								\abs{z_{u_m}(\lam_{u_m+1})}&= \mu\abs{\sum_{\iota = \ell_m}^{u_m} {\ip{\f_{\iota}(\lam_{u_m+1}) - \f(\lam_{u_m+1}),\lam_{u_m+1}-\ls}}} \\
								&\hspace{-1cm}\leq  \norm{\mu\sum_{\iota = \ell_m}^{u_m} \left(\f_{\iota}(\lam_{u_m+1}) - \f(\lam_{u_m+1})\right)}\norm{\lam_{u_m+1}-\ls} \leq L^1_{1/\mu}(\T_m)\Lambda_{\max}. \label{slbnd}
								\end{align}
								Similarly, denoting $\T'_\tau := \{\ell_m,\ell_m+1, \ldots, \tau\}$ for all $\ell_m \leq \tau \leq u_m$, it holds from using triangle inequality and \eqref{a5a6}, that
								\small
								\begin{align}
								&\!\!\!\abs{z_\tau(\lam_{\tau+1})-z_\tau(\lam_\tau)} \!=\! \epsilon\Biggr|\sum_{\iota  = \ell_m}^\tau  {\ip{\f_{\iota}(\lam_{\tau + 1}) \!- \!\f(\lam_{\tau+1}),\lam_{\tau+1}\!\!-\!\!\ls}}-  {\ip{\f_{\iota}(\lam_{\tau}) - \f(\lam_{\tau}),\lam_{\tau}-\ls}}\Biggr| 
								\\
								&\hspace{-.2cm} \leq \epsilon \Biggr|\sum_{\iota  = \ell_m}^\tau  \!\!\langle\f_{\iota}(\lam_{\tau + 1}\!) \!-\! \f(\lam_{\tau+1})\! -\! \f_{\iota}(\lam_{\tau}) \!+\! \f(\lam_{\tau}), (\lam_{\tau+1}\!-\!\ls)\rangle\Biggr|  + \mu \abs{\sum_{\iota = \ell_m}^{\tau}{\ip{\f_{\iota}(\lam_{\tau}) - \f(\lam_{\tau}),\lam_{\tau+1} - \lam_\tau}}} \nonumber
								\\
								& \hspace{-0.2cm}\leq  \norm{\epsilon\!\!\sum_{\iota  = \ell_m}^\tau \!\! \f_{\iota}(\lam_{\tau + 1})\! -\! \f(\lam_{\tau+1}) \!-\! {\f_{\iota}(\lam_{\tau})\!+\! \f(\lam_{\tau})}} \norm{\lam_{\tau+1}-\ls} {+} \norm{\epsilon\sum_{\iota = \ell_m}^{\tau}\left(\f_{\iota}(\lam_{\tau}) - \f(\lam_{\tau})\right)}\norm{\lam_{\tau+1} - \lam_\tau}\nonumber
								\\
								&\hspace{-0.2cm} \leq \left(L^2_{1/\epsilon}(\T'_\tau)\Lambda_{\max} + L^1_{1/\epsilon}(\T'_\tau)\right)\norm{\lam_{\tau+1}-\lam_\tau}  \leq \epsilon\left(L^2_{1/\epsilon}(\T'_\tau)\Lambda_{\max} + L^1_{1/\epsilon}(\T'_\tau)\right)G \label{slineq} 
								\end{align}
								\normalsize
								where the \eqref{slineq} uses the non-expansive property of the projection operator $\pl{\cdot}$ and the boundedness of the stochastic subgradients (cf. (\textbf{A3}$^\prime$)). Substituting \eqref{slbnd} and \eqref{slineq} into the expression for $C^m(\epsilon)$ yields the following bound
								\small
								\begin{align}
								\abs{C^m(\epsilon)} &\leq L^1_{1/\mu}(\T_m)\Lambda_{\max} + \epsilon G\sum_{\tau = \ell_m}^{u_m} L^2_{1/\epsilon}(\T'_\tau)\Lambda_{\max} + L^1_{1/\epsilon}(\T'_\tau) \nonumber\\
								&\leq L^1_{1/\mu}(\T_m)\Lambda_{\max}  + G\sup_{\ell_m \leq \tau \leq u_m} (L^2_{1/\epsilon}(\T'_\tau)\Lambda_{\max} + L^1_{1/\epsilon}(\T'_\tau)) \nonumber.
								\end{align}
								\normalsize
								Finally, the bound for $C_t(n,\epsilon)$ becomes
								\small
								\begin{align}
								C_t(n,\epsilon) &\leq \frac{1}{n}\sum_{m=0}^{n} \abs{C^m(\epsilon)} \leq  \sup_{0\leq  m \leq n} \abs{C^m(\epsilon)}  \\
								&\leq \sup_{0\leq m \leq n} L^1_{1/\epsilon}(\T_m)\Lambda_{\max}+ G\sup_{0\leq \tau < t} (L^2_{1/\epsilon}(\T'_\tau)\Lambda_{\max} + L^1_{1/\epsilon}(\T'_\tau)) \label{ctbound}
								\end{align}		
								\normalsize
								Therefore, the rate result from Lemma \ref{errbounds} implies that 
								\small
								\begin{align}
								 \epsilon^{\zeta-1/2}C_t(n,\epsilon) &\leq \Lambda_{\max}\sup_{0\leq m \leq n} \epsilon^{\zeta-1/2}L^1_{1/\epsilon}(\T_m)+ {G}\sup_{0\leq \tau < t} (\epsilon^{\zeta-1/2}L^2_{1/\epsilon}(\T'_\tau)\Lambda_{\max} + \epsilon^{\zeta-1/2}L^1_{1/\epsilon}(\T'_\tau))
								\end{align}	
								\normalsize
								which goes to zero almost surely as $\epsilon \rightarrow 0$, yielding the bound in \eqref{ctlimep}. {Likewise, let $A_{1m} < \infty$ be the constant associated with the bounds for $\Lambda_{\max}L_{1/\epsilon}^1(\T_m)$, as necessitated by Lemma \ref{errbounds}. Then, using the union bound, it follows that
									\begin{align}
									\mathbb{P}\left(\abs{\sup_{0\leq m \leq n} L_{1/\epsilon}^1(\T_m)} > \epsilon^{1/2-\zeta}\right) &\leq \sum_{m=0}^n A_{1m}\exp\left(-\epsilon^{-2\zeta}\right) \leq A_1\exp\left(-\epsilon^{-2\zeta}\right)
									\end{align}
									where $A_1 := \sum_{m} A_{1m}$. Along the same lines, the result in Lemma \ref{errbounds} and the subsequent use of the union bound imply that there exist a constant $A_2 < \infty$ such that the probability of the second term in \eqref{ctbound} exceeds $\epsilon^{1/2-\zeta}$ is bounded by $A_1\exp\left(-\epsilon^{-2\zeta}\right)$. Combining the two bounds, and again using union bound, we have that 
									\begin{align}\label{first_prob_bound}
									\mathbb{P}\left(C_{t}(n,\epsilon) > \epsilon^{1/2-\zeta}\right) \leq A_3\exp\left(-\epsilon^{-2\zeta}\right)
									\end{align} 
									where $A_3 = A_1+A_2$.}
								
								\emph{Fixed $\epsilon > 0$ and $n \rightarrow \infty$: } 
								In this case, $C_t(n,\epsilon)$ must now be split into two terms as follows,
								\begin{align}
								{{C_t(n,\epsilon)}} \leq \epsilon \abs{C(n)} + \frac{{\abs{D(n)}}}{n}
								\end{align}
								where,\vspace{-6mm}								
								\small
								\begin{align}
								\hspace{1cm}C(n) &:=  \sum_{\tau = \ell_m}^{n/\mu-1}{\ip{\f_\tau(\lam_\tau)-\f(\lam_\tau),\lam_\tau-\ls}} \\
								&\hspace{-22mm}= \sum_{\tau = 0}^{1/\mu-1}\frac{1}{n} \sum_{m=0}^{n-1} {\ip{\f_{m/\mu + \tau}(\lam_{m/\mu + \tau})-\f(\lam_{m/\mu+\tau}),\lam_{m/\mu+\tau}-\ls}} \nonumber\\
								\hspace{-5mm}D(n) &:= \mu \sum_{\tau = n/\mu}^{t-1} {\ip{\f_\tau(\lam_\tau)-\f(\lam_\tau),\lam_\tau-\ls}}.
								\end{align}
								\normalsize
								For this analysis, it is assumed without loss of generality that $1/\epsilon$ is an integer. That way, the subscripts $\frac{m}{\epsilon}+ \tau$ are also integers and the floor operation is not required. Given $\epsilon$, note that $D(n)$ is a sum of a fixed number of bounded random variables, so that $D(n)/n \rightarrow 0$ surely as $n \rightarrow \infty$.
								In order to bound $C(n)$, define for all $\lam \in\Lambda$ and $0  \leq \tau \leq 1/\mu-1$, 
								\begin{align}
								z_{\tau}(\lam) = \sum_{\iota = 0}^{\tau} \frac{1}{n}\sum_{m=0}^{n-1}{\ip{\f_{m/\epsilon+\tau}(\lam) - \f(\lam),\lam-\ls}}.
								\end{align}
								Then, it follows that
								\begin{align}
								C(n) = z_{n/\mu-1}(\lam_{n/\mu}) - \sum_{\tau = 0}^{1/\mu-1} \left(z_\tau(\lam_{\tau+1}) - z_\tau(\lam_\tau)\right). \label{cnsplit}
								\end{align}
								It is now possible to bound each term in \eqref{cnsplit} separately. Defining $\T^{\tau}:=\{m/\mu + \tau\}_{m=0}^{n-1}$, and using \eqref{a5a6}, it follows that
								\begin{align}
								\abs{z_{n/\mu-1}(\lam_{n/\mu})}&= \sum_{\tau = 0}^{1/\mu-1} \abs{\frac{1}{n}\sum_{m=0}^{n-1}{\ip{\f_{m/\mu+\tau}(\lam_{n/\mu}) - \f(\lam_{n/\mu}),\lam_{n/\mu}-\ls}}}\nonumber\\
								&\leq  \sum_{\tau = 0}^{1/\mu-1}\norm{\frac{1}{n}\sum_{m=0}^{n-1} \left(\f_{m/\mu+\tau}(\lam_{n/\mu}) - \f(\lam_{n/\mu})\right)}\norm{\lam_{n/\mu}-\ls} \nonumber\\
								&\leq \Lambda_{\max} \sum_{\tau = 0}^{1/\mu-1}{L^1_{n}}(\T^\tau). \label{slbnd2}
								\end{align}
								Proceeding similarly,
								\small
								\begin{align}
								\abs{z_\tau(\lam_{\tau+1})-z_\tau(\lam_\tau)}   &=\sum_{\iota = 0}^{\tau} \Biggr|\frac{1}{n}\sum_{m  = 0}^{n-1}  {\ip{\f_{m/\mu+\iota}(\lam_{\tau + 1}) - \f(\lam_{\tau+1}),\lam_{\tau+1}-\ls}} -  {\ip{\f_{m/\mu+\iota}(\lam_{\tau}) - \f(\lam_{\tau}),\lam_{\tau}-\ls}}\Biggr|\nonumber\\	
								& \leq  \sum_{\iota = 0}^{\tau} \Biggr|\Biggr|\frac{1}{n}\sum_{m  = 0}^{n-1}  (\f_{m/\mu+\iota}(\lam_{\tau + 1}) - \f(\lam_{\tau+1})  - \f_{m/\mu+\iota}(\lam_{\tau }) + \f(\lam_{\tau}))\Biggr|\Biggr|\norm{\lam_{\tau+1}-\ls} \nonumber\\
								& \quad+\sum_{\iota = 0}^{\tau}\norm{\frac{1}{n}\sum_{m = 0}^{n-1}\left(\f_{m/\mu+\iota}(\lam_{\tau}) - \f(\lam_{\tau})\right)}\norm{\lam_{\tau+1} - \lam_\tau}\nonumber\\
								&\leq \sum_{\iota = 0}^\tau\left(L^2_{n}(\T^\iota)\Lambda_{\max} + L^1_{n}(\T^\iota)\right)\norm{\lam_{\tau+1}-\lam_\tau} \label{sl2}\\
								&\leq \epsilon G \sum_{\iota = 0}^\tau\left(L^2_{n}(\T^\iota)\Lambda_{\max} + L^1_{n}(\T^\iota)\right). \label{slineq2}
								\end{align}
								\normalsize
								Finally,  substituting \eqref{slbnd2} and \eqref{slineq2} into the expression for $C(n)$, and noting that $1/\epsilon$ is a fixed number, the following bound is obtained
								\begin{align}
								\epsilon\abs{C(n)} &\leq \Lambda_{\max} \epsilon\sum_{\tau = 0}^{1/\mu-1}L^1_{n}(\T^\tau) + \epsilon^2\sum_{\tau = 0}^{1/\epsilon-1} G \sum_{\iota = 0}^\tau\big(L^2_{n}(\T^\iota)\Lambda_{\max} + L^1_{n}(\T^\iota)\big)\nonumber
								\end{align}
								\begin{align}
								&\leq \Lambda_{\max} \sup_{0\leq \tau < 1/\epsilon} L^2_{n}(\T^\tau)+ G\sup_{0\leq \tau < 1/\epsilon}\sup_{0\leq \iota \leq \tau} \big(L^2_{n}(\T^\iota)\Lambda_{\max}  + L^1_{n}(\T^\iota)\big) \label{cnbnd}
								\end{align}
								which goes to zero almost surely as $n \rightarrow \infty$, implying that $C_t(n,\epsilon) \rightarrow 0$ almost surely as $n \rightarrow \infty$. Both the rate results can again be inferred as in the previous case. Indeed, similar to \eqref{first_prob_bound}, given $\zeta>0$, there exist $A_4<\infty$ such that $\mathbb{P}\left(C_{t}(n,\epsilon) > n^{\zeta-1/2}\right) \leq A_4\exp(-n^{2\zeta})$. Combining with \eqref{first_prob_bound}, the probability bounds can be written as
								\begin{align}\label{probbounds2}
								\mathbb{P}\left(C_{t}(n,\epsilon) > \epsilon^{1/2-\zeta}\right) \leq A\exp(-\epsilon^{-2\zeta})\ \ \text{and}\ \ 	\mathbb{P}\left(C_{t}(n,\epsilon) > n^{\zeta-1/2}\right) \leq A\exp(-n^{2\zeta})
								\end{align}
								where $A\!\! = \!\max\{A_3,\!A_4\}$. The required result follows by choosing $\nu = \max\{\frac{1}{\epsilon},n\}$ in \eqref{probbounds2}.  
							\end{IEEEproof}
							
							\vspace{-5mm}
							\section{Derivation of \eqref{itlam1} and \eqref{itlam2}} \label{itlam}
							Consider first the slow fading case, where the winning user is given by
							\vspace{-5mm}
							\begin{align}
							i_t&=\argmax_{{i\in\mathcal{M}_t}}  \lambda_t \left(W\colb{\log_2}(\hat{p}^i_t\gamma^i_t/\alpha)\right)-c_t^i\hat{p}^i_t	\label{itlams1}
							\end{align}
							where $\hat{p}^i_t$ is given by \eqref{primal_slow2}. Thus, the objective function in \eqref{itlams1} for a given $\lambda$ can be written as
							\small
							\begin{align}\label{condition}
							T^i_t(\lambda)=\begin{cases}
							\lambda W\colb{\log_2}(C_{\min}\frac{\gamma^i_t}{c^i_t})-C_{\min}, & {\lambda\leq\frac{C_{\min}}{W}} \\
							\lambda W\colb{\log_2}(C_{\max}\frac{\gamma^i_t}{c^i_t})-C_{\max}, & \text{if}\ \lambda \geq \frac{C_{\max}}{W} \\
							\lambda W\colb{\log_2}(\frac{\lambda W}{\alpha}\frac{\gamma^i_t}{c^i_t})-\lambda W, & \text{otherwise}. 
							\end{cases}
							\end{align}
							\normalsize
							Since $\colb{\log_2}$ is monotonic function, observe in \eqref{condition} that in all three cases, $T^i_t(\lambda)$ depends monotonically on $\gamma^i_t/c^i_t$ for all $\lambda > 0$. This allows us to conclude that $i_t = \argmax_{{i\in\mathcal{M}_t}} T^i_t(\lambda)  = \argmax_{{i\in\mathcal{M}_t}} \frac{\gamma^i_t}{c^i_t}$	which is the required identity in \eqref{itlam1}.							
							Similarly for the fast fading case, the objective function for the winning user in \eqref{itlam2} is given by
							\small
							\begin{align}\label{condition2}
							T^i_t(\lambda)\!=\!\!\begin{cases}
							\lambda W\colb{\log_2}(C_{\min}\frac{\gamma^i_t}{c^i_t})\!\! +\!\! \lambda W\psi_i\!-\!C_{\min}, & \text{if}\ {\lambda\leq \frac{C_{\min}}{W}} \\
							\lambda W\colb{\log_2}(C_{\max}\frac{\gamma^i_t}{c^i_t}) + \lambda W\psi_i-C_{\max}, & \text{if}\ \lambda \geq \frac{C_{\max}}{W} \\
							\lambda W\colb{\log_2}(\frac{\lambda W}{\alpha}\frac{\gamma^i_t}{c^i_t}) + \lambda W\psi_i-\lambda W, & \hspace{0mm}\!\!\!\!\!\!\!\!\!\!\! \text{otherwise}
							\end{cases}
							\end{align}\normalsize
							which, for $\lambda > 0$, again depends monotonically on $\colb{\log_2}(\gamma^i_t/c^i_t) + \psi_i$. The expression in \eqref{itlam2} therefore follows.
							\vspace{-8mm}
							\footnotesize
							\bibliographystyle{IEEEtran} 
							\bibliography{IEEEabrv,reference}

\begin{thebibliography}{10}
\providecommand{\url}[1]{#1}
\csname url@samestyle\endcsname
\providecommand{\newblock}{\relax}
\providecommand{\bibinfo}[2]{#2}
\providecommand{\BIBentrySTDinterwordspacing}{\spaceskip=0pt\relax}
\providecommand{\BIBentryALTinterwordstretchfactor}{4}
\providecommand{\BIBentryALTinterwordspacing}{\spaceskip=\fontdimen2\font plus
\BIBentryALTinterwordstretchfactor\fontdimen3\font minus
  \fontdimen4\font\relax}
\providecommand{\BIBforeignlanguage}[2]{{%
\expandafter\ifx\csname l@#1\endcsname\relax
\typeout{** WARNING: IEEEtran.bst: No hyphenation pattern has been}%
\typeout{** loaded for the language `#1'. Using the pattern for}%
\typeout{** the default language instead.}%
\else
\language=\csname l@#1\endcsname
\fi
#2}}
\providecommand{\BIBdecl}{\relax}
\BIBdecl

\bibitem{resource_alloca_tut_1_neely}
L.~Georgiadis, M.~J. Neely, and L.~Tassiulas, ``Resource allocation and
  cross-layer control in wireless networks,'' \emph{Found. Trends. Network.},
  vol.~1, no.~1, pp. 1--149, 2006.

\bibitem{alloc_smart_tut_1}
Z.~Fan, P.~Kulkarni, S.~Gormus, C.~Efthymiou, G.~Kalogridis, M.~Sooriyabandara,
  Z.~Zhu, S.~Lambotharan, and W.~H. Chin, ``Smart grid communications: overview
  of research challenges, solutions, and standardization activities,''
  \emph{{IEEE} Commun. Surveys Tuts.}, vol.~15, no.~1, pp. 21--38, 2013.

\bibitem{alloc_sched_tut_1}
R.~Afolabi, A.~Dadlani, and K.~Kim, ``{M}ulticast {S}cheduling and {R}esource
  {A}llocation {A}lgorithms for {OFDMA}-{B}ased {S}ystems: {A} {S}urvey,''
  \emph{{IEEE} Commun. Surveys Tuts.}, vol.~15, no.~1, pp. 240--254, First
  2013.

\bibitem{stoch_res_prob_3}
X.~Wang and N.~Gao, ``Stochastic resource allocation over fading multiple
  access and broadcast channels,'' \emph{IEEE Trans. Inf. Theory}, vol.~56,
  no.~5, pp. 2382--2391, 2010.

\bibitem{separation}
A.~Ribeiro and G.~Giannakis, ``Separation {P}rinciples in {W}ireless
  {N}etworking,'' \emph{IEEE Trans. Inf. Theory}, vol.~56, no.~9, pp.
  4488--4505, Sept. 2010.

\bibitem{marques3}
A.~G. Marques, L.~M. Lopez-Ramos, G.~B. Giannakis, and J.~Ramos, ``Resource
  allocation for interweave and underlay {CR}s under
  probability-of-interference constraints,'' \emph{{IEEE} J. Sel. Areas
  Commun.}, vol.~30, no.~10, pp. 1922--1933, 2012.

\bibitem{crolayketan}
K.~Rajawat, N.~Gatsis, and G.~Giannakis, ``Cross-{L}ayer {D}esigns in {C}oded
  {W}ireless {F}ading {N}etworks {W}ith {M}ulticast,'' \emph{IEEE/ACM Trans.
  Netw.}, vol.~19, no.~5, pp. 1276--1289, Oct 2011.

\bibitem{ale10}
A.~Ribeiro, ``{E}rgodic {S}tochastic {O}ptimization {A}lgorithms for {W}ireless
  {C}ommunication and {N}etworking,'' \emph{{IEEE} Trans. Signal Process.},
  vol.~58, no.~12, pp. 6369--6386, Dec. 2010.

\bibitem{ofdm}
X.~Wang and G.~Giannakis, ``Resource {A}llocation for {W}ireless {M}ultiuser
  {OFDM} {N}etworks,'' \emph{IEEE Trans. Inf. Theory}, vol.~57, no.~7, pp.
  4359--4372, July 2011.

\bibitem{nedic2001convergence}
A.~Nedi{\'c} and D.~Bertsekas, ``Convergence rate of {I}ncremental
  {S}ubgradient {A}lgorithms,'' in \emph{Stochastic {O}ptimization:
  {A}lgorithms and {A}pplications}.\hskip 1em plus 0.5em minus 0.4em\relax
  Springer, 2001, pp. 223--264.

\bibitem{bertsekas2011incremental}
D.~P. Bertsekas, ``Incremental {G}radient, {S}ubgradient, and {P}roximal
  methods for {C}onvex {O}ptimization: {A} survey,'' \emph{LIDS-P-2848}, 2011.

\bibitem{Luo_LMS_rate}
Z.-Q. Luo, ``On the {C}onvergence of the {LMS} {A}lgorithm with {A}daptive
  {L}earning {R}ate for {L}inear {F}eedforward {N}etworks,'' \emph{Neural
  Comput.}, vol.~3, no.~2, pp. 226--245, Jun. 1991.

\bibitem{bottou1998online}
L.~Bottou, ``Online learning and stochastic approximations,'' \emph{On-line
  learning in neural networks}, vol.~17, no.~9.

\bibitem{bottou2012stochastic}
------, ``Stochastic gradient descent tricks,'' in \emph{Neural Networks:
  Tricks of the Trade}.\hskip 1em plus 0.5em minus 0.4em\relax Springer, 2012,
  pp. 421--436.

\bibitem{bach2011non}
E.~Moulines and F.~R. Bach, ``Non-asymptotic analysis of stochastic
  approximation algorithms for machine learning,'' in \emph{Adv. Neural Inf.
  Process. Syst.}, 2011, pp. 451--459.

\bibitem{nedic2009approximate}
A.~Nedic and A.~Ozdaglar, ``Approximate primal solutions and rate analysis for
  dual subgradient methods,'' \emph{SIAM J. Optim.}, vol.~19, no.~4, pp.
  1757--1780, 2009.

\bibitem{ji2016wireless}
M.~Ji, G.~Caire, and A.~F. Molisch, ``Wireless device-to-device caching
  networks: Basic principles and system performance,'' \emph{IEEE J. Sel. Areas
  Commun}, vol.~34, no.~1, pp. 176--189, 2016.

\bibitem{archi_D2D}
N.~Golrezaei, P.~Mansourifard, A.~F. Molisch, and A.~G. Dimakis, ``Base-station
  assisted device-to-device communications for high-throughput wireless video
  networks,'' \emph{IEEE Trans. Wireless Commun.}, vol.~13, no.~7, pp.
  3665--3676, 2014.

\bibitem{robbins1951}
H.~Robbins and S.~Monro, ``A {S}tochastic {A}pproximation {M}ethod,''
  \emph{Ann. Math. Statistics}, vol.~22, no.~3, pp. 400--407, 09 1951.

\bibitem{Widrow}
B.~Widrow and S.~D. Stearns, \emph{Adaptive {S}ignal {P}rocessing}.\hskip 1em
  plus 0.5em minus 0.4em\relax Upper Saddle River, NJ, USA: Prentice-Hall,
  Inc., 1985.

\bibitem{sayed2011adaptive}
A.~H. Sayed, \emph{Adaptive filters}.\hskip 1em plus 0.5em minus 0.4em\relax
  John Wiley \& Sons, 2011.

\bibitem{bottou_91c_neural_learning}
L.~Bottou, ``{S}tochastic {G}radient {L}earning in {N}eural {N}etworks,'' in
  \emph{Proceedings of Neuro-N\^imes'91}.\hskip 1em plus 0.5em minus
  0.4em\relax Nimes, France: EC2, 1991.

\bibitem{kushner1994analysis}
H.~J. Kushner and J.~Yang, ``Analysis of adaptive step size {SA} algorithms for
  parameter tracking,'' in \emph{Proc. of the IEEE Conf. on Decision and
  Control}, vol.~1, 1994, pp. 730--737.

\bibitem{bottou2010large}
L.~Bottou, ``Large-scale machine learning with stochastic gradient descent,''
  in \emph{Proc. of COMPSTAT}.\hskip 1em plus 0.5em minus 0.4em\relax Springer,
  2010, pp. 177--186.

\bibitem{optml}
S.~Suvrit, S.~Nowozin, and S.~J. Wright, Eds., \emph{Optimization for Machine
  Learning}.\hskip 1em plus 0.5em minus 0.4em\relax Cambirdge, Mass.: MIT
  Press, 2012.

\bibitem{kelly1997charging}
F.~Kelly, ``Charging and rate control for elastic traffic,'' \emph{Euro. Trans.
  Telecommun.}, vol.~8, no.~1, pp. 33--37, 1997.

\bibitem{larsson1999ergodic}
T.~Larsson, M.~Patriksson, and A.-B. Str{\"o}mberg, ``Ergodic, primal
  convergence in dual subgradient schemes for convex programming,'' \emph{Math.
  Program.}, vol.~86, no.~2, pp. 283--312, 1999.

\bibitem{classconv}
N.~Gatsis, A.~Ribeiro, and G.~Giannakis, ``A class of convergent algorithms for
  resource allocation in wireless fading networks,'' \emph{IEEE Trans. Wireless
  Commun.}, vol.~9, no.~5, pp. 1808--1823, May 2010.

\bibitem{Solo_1994}
V.~Solo and X.~Kong, \emph{{A}daptive {S}ignal {P}rocessing {A}lgorithms:
  {S}tability and {P}erformance}.\hskip 1em plus 0.5em minus 0.4em\relax Upper
  Saddle River, NJ, USA: Prentice-Hall, 1995.

\bibitem{neely2010stochastic}
M.~J. Neely, ``{S}tochastic network optimization with application to
  communication and queueing systems,'' \emph{Synthesis Lectures on Commun.
  Netw.}, vol.~3, no.~1, pp. 1--211, 2010.

\bibitem{yu2015convergence}
H.~Yu and M.~J. Neely, ``On the convergence time of the drift-plus-penalty
  algorithm for strongly convex programs,'' in \emph{Proc. of the IEEE Conf. on
  Decision and Control}, 2015, pp. 2673--2679.

\bibitem{sidiropoulos2006transmit}
N.~D. Sidiropoulos, T.~N. Davidson, and Z.-Q.~T. Luo, ``Transmit beamforming
  for physical-layer multicasting,'' \emph{{IEEE} Trans. Signal Process.},
  vol.~54, no.~6, pp. 2239--2251, 2006.

\bibitem{resource_alloca_tut_4}
D.~Palomar and M.~Chiang, ``A tutorial on decomposition methods for network
  utility maximization,'' \emph{{IEEE} J. Sel. Areas Commun.}, vol.~24, no.~8,
  pp. 1439--1451, Aug 2006.

\bibitem{gatsis2014stochastic}
N.~Gatsis and A.~G. Marques, ``A stochastic approximation approach to load
  shedding in power networks,'' in \emph{Proc. of the IEEE ICASSP}, 2014, pp.
  6464--6468.

\bibitem{wang2016dynamic}
X.~Wang, Y.~Zhang, T.~Chen, and G.~B. Giannakis, ``Dynamic energy management
  for smart-grid-powered coordinated multipoint systems,'' \emph{IEEE J. Sel.
  Areas Commun.}, vol.~34, no.~5, pp. 1348--1359, 2016.

\bibitem{wang2016dynamic2}
X.~Wang, T.~Chen, X.~Chen, X.~Zhou, and G.~B. Giannakis, ``Dynamic resource
  allocation for smart-grid powered mimo downlink transmissions,'' \emph{IEEE
  J. Sel. Areas Commun.}, vol.~34, no.~12, pp. 3354--3365, 2016.

\bibitem{wang2016stochastic}
X.~Wang, T.~Ma, R.~Zhang, and X.~Zhou, ``Stochastic online control for
  energy-harvesting wireless networks with battery imperfections,'' \emph{IEEE
  Trans. Wireless Commun.}, vol.~15, no.~12, pp. 8437--8448, 2016.

\bibitem{tse}
D.~Tse and P.~Viswanath, \emph{Fundamentals of wireless communication}.\hskip
  1em plus 0.5em minus 0.4em\relax Cambridge University Press, 2005.

\bibitem{Gatsis_2}
N.~Gatsis and G.~Giannakis, ``Residential {L}oad {C}ontrol: {D}istributed
  {S}cheduling and {C}onvergence {W}ith {L}ost {AMI} {M}essages,'' \emph{IEEE
  Trans. Smart Grid}, vol.~3, no.~2, pp. 770--786, June 2012.

\bibitem{rao1972remark}
K.~B. Rao, M.~B. Rao \emph{et~al.}, ``A remark on nonatomic measures,''
  \emph{Ann. Math. Stat.}, vol.~43, no.~1, pp. 369--370, 1972.

\bibitem{hermes1969functional}
H.~Hermes and J.~P. Lasalle, \emph{Functional analysis and time optimal
  control}.\hskip 1em plus 0.5em minus 0.4em\relax Academic Press, 1969.

\bibitem{shitovitz1973oligopoly}
B.~Shitovitz, ``Oligopoly in markets with a continuum of traders,''
  \emph{Econometrica}, pp. 467--501, 1973.

\bibitem{luo2008dynamic}
Z.-Q. Luo and S.~Zhang, ``Dynamic spectrum management: Complexity and
  duality,'' \emph{IEEE J. Sel. Topics Signal Process.}, vol.~2, no.~1, pp.
  57--73, 2008.

\bibitem{toh2010accelerated}
K.-C. Toh and S.~Yun, ``An accelerated proximal gradient algorithm for nuclear
  norm regularized linear least squares problems,'' \emph{Pacific J.Opt.},
  vol.~6, no. 615-640, p.~15, 2010.

\bibitem{recht2011hogwild}
B.~Recht, C.~Re, S.~Wright, and F.~Niu, ``Hogwild: A lock-free approach to
  parallelizing stochastic gradient descent,'' in \emph{Advances in Neural
  Information Processing Systems}, 2011, pp. 693--701.

\bibitem{slln}
V.~V. Petrov, \emph{Limit {T}heorems of {P}robability {T}heory}.\hskip 1em plus
  0.5em minus 0.4em\relax Oxford Studies in Probability, 1995.

\bibitem{fuk1971probability}
D.~K. Fuk and S.~V. Nagaev, ``Probability inequalities for sums of independent
  random variables,'' \emph{Theory of Probability \& Its Applications},
  vol.~16, no.~4, pp. 643--660, 1971.

\bibitem{baum1962exponential}
L.~E. Baum, M.~Katz, and R.~R. Read, ``Exponential convergence rates for the
  law of large numbers,'' \emph{Trans. Amer. Math. Soc.}, vol. 102, no.~2, pp.
  187--199, 1962.

\end{thebibliography}

						\end{document}